\documentclass[graybox]{svmult}

\usepackage{type1cm}
\usepackage{makeidx}           
\usepackage{graphicx}          
\usepackage{multicol}          
\usepackage[bottom]{footmisc} 
\usepackage[T1]{fontenc}

\usepackage{hyperref}
\hypersetup{colorlinks=true, urlcolor=blue, citecolor=cyan}
\usepackage[comma,authoryear]{natbib}
\usepackage[english]{babel}
\usepackage{amsmath}
\usepackage{amssymb}
\usepackage{booktabs}

\usepackage{amsfonts}
\usepackage{caption}
\usepackage{subcaption}
\usepackage{suffix}
\usepackage{mathtools}
\usepackage{paralist}
\usepackage{bigints}
\usepackage{bbm}

\DeclarePairedDelimiterX\MeijerM[3]{\lparen}{\rparen}%
{\begin{smallmatrix}#1 \\ #2\end{smallmatrix}\delimsize\vert\,#3}

\newcommand\MeijerG[8][]{%
	G^{\,#2,#3}_{#4,#5}\MeijerM[#1]{#6}{#7}{#8}}
\WithSuffix\newcommand\MeijerG*[7]{%
	G^{\,#1,#2}_{#3,#4}\MeijerM*{#5}{#6}{#7}}
\long\def\symbolfootnote[#1]#2{\begingroup%
	\def\thefootnote{\fnsymbol{footnote}}\footnote[#1]{#2}\endgroup}

\DeclareMathOperator{\EV}{\mathbb{E}}

\DeclareMathOperator{\Ei}{Ei}
\DeclareMathOperator{\E1}{E_1}

\DeclareMathOperator{\SADD}{SADD}

\newcommand{\iu}{\mathrm{i}\mkern1mu}
\renewcommand{\le}{\leqslant} 
\renewcommand{\ge}{\geqslant}
\newcommand{\abs}[1]{\left\vert#1\right\vert}
\DeclareMathOperator{\One}{\mathbbm{1}}
\newcommand{\indicator}[1]{{\One_{\left\{#1\right\}}}}

\makeindex

\begin{document}
	
\bibliographystyle{spbasic}

\title*{More On the Quasi-Stationary Distribution of the Shiryaev--Roberts Diffusion}

\titlerunning{Quasi-Stationary Distribution of the SR Diffusion}

\author{Soumik Banerjee and Aleksey S. Polunchenko}

\institute{
Soumik Banerjee
\at Department of Mathematics and Statistics, State University of New York at Binghamton, USA, \email{soumik@math.binghamton.edu}
\and
Aleksey S. Polunchenko
\at Department of Mathematics and Statistics, State University of New York at Binghamton, USA,	\email{aleksey@binghamton.edu}
}

\maketitle

\abstract{ We consider the classical Shiryaev--Roberts martingale diffusion, $(R_t)_{t\ge0}$, restricted to the interval $[0,A]$, where $A>0$ is a preset absorbing boundary. We take yet another look at the well-known phenomenon of quasi-stationarity (time-invariant probabilistic behavior, conditional on no absorbtion hitherto) exhibited by the diffusion in the temporal limit, as $t\to+\infty$, for each $A>0$. We obtain new upper- and lower-bounds for the quasi-stationary distribution's probability density function (pdf), $q_{A}(x)$; the bounds vary in the trade-off between simplicity and tightness. The bounds imply directly the expected result that $q_{A}(x)$ converges to the pdf, $h(x)$, of the diffusion's stationary distribution, as $A\to+\infty$; the convergence is pointwise, for all $x\ge0$. The bounds also yield an explicit upperbound for the gap between $q_{A}(x)$ and $h(x)$ for a fixed $x$. By virtue of integration the bounds for the pdf $q_{A}(x)$ translate into new bounds for the corresponding cumulative distribution function (cdf), $Q_{A}(x)$. All of our results are established explicitly, using certain latest monotonicity properties of the modified Bessel $K$ function involved in the exact closed-form formula for $q_{A}(x)$ recently obtained by~\cite{Polunchenko:SA2017a}. We conclude with a discussion of potential applications of our results in quickest change-point detection: our bounds allow for a very accurate performance analysis of the so-called randomized Shiryaev--Roberts--Pollak change-point detection procedure.
}

\keywords{
Generalized Shiryaev--Roberts procedure;
Markov diffusion;
Quasi-stationary distribution;
Quickest change-point detection;
Whittaker functions;
Modified Bessel functions.}


\section{Introduction} 

This work is a continuation of the recent paper by~\cite{Li+Polunchenko:SA2020} and, too, focuses on the phenomena of quasi-stationarity and stationarity exhibited by one particular version of the Generalized Shiryaev--Roberts (GSR) stochastic process---a time-homogeneous Markov diffusion well-known in the area of quickest change-point detection. See, e.g.,~\cite{Shiryaev:SMD61,Shiryaev:TPA63,Shiryaev:Book78,Shiryaev:Bachelier2002,Shiryaev:Book2011,Shiryaev:Book2017},~\cite{Pollak+Siegmund:B85},~\cite{Feinberg+Shiryaev:SD2006},~\cite{Burnaev+etal:TPA2009},~\cite{Polunchenko+Sokolov:MCAP2016}, and~\cite{Polunchenko:SA2016,Polunchenko:SA2017a,Polunchenko:SA2017b,Polunchenko:TPA2017}. More specifically, the GSR process' version of interest is the solution $(R_{t}^{r})_{t\ge0}$ of the stochastic differential equation
\begin{align} \label{eq:Rt_r-def}
  dR_{t}^{r} & = dt+R_{t}^{r} dB_{t} \;\text{with}\; R_{0}^{r} \coloneqq r\ge0
\;
\text{fixed},
\end{align}
where $(B_{t})_{t\ge0}$ is standard Brownian motion in the sense that $\EV[dB_t]=0$, $\EV[(dB_t)^2]=dt$, and $B_0=0$; the initial value $R_{0}^{r}\coloneqq r\ge0$ is often referred to as the process' headstart. It is straightforward to solve~\eqref{eq:Rt_r-def} and express $(R_{t}^{r})_{t\ge0}$ explicitly as
\begin{align*}
R_{t}^{r}
&=
\exp\left\{B_{t}-\dfrac{1}{2}t\right\}\left(r+\bigintsss_{0}^{t} \exp\biggl\{-\left(B_s-\dfrac{1}{2}s\right)\biggr\}ds\right),
\;\;
t\ge0,
\end{align*}
so that the set $[0,+\infty)$ is easily seen to be the ``natural'' state space for $(R_{t}^{r})_{t\ge0}$ because $R_0^{r}\coloneqq r\ge0$ by assumption. Moreover, it is also easily checked that $\EV[R_{t}^{r}-t-r]=0$ for any $t,r\ge0$, i.e., the process $\{R_{t}^{r}-t-r\}_{t\ge0}$ is a zero-mean martingale. Yet, although $(R_{t}^{r})_{t\ge0}$ has a linear upward trend in time, it is actually a recurrent process with a nontrivial probabilistic behavior in the limit as $t\to+\infty$; cf.~\cite[p.~270]{Pollak+Siegmund:B85}. Specifically, if $(R_{t}^{r})_{t\ge0}$ is let run ``loose'', i.e., considered on the entire nonnegative half-line, then the limiting (as $t\to+\infty$) behavior of $(R_{t}^{r})_{t\ge0}$ is known as stationarity, and it is characterized by the invariant probability measure whose cumulative distribution function (cdf) and density (pdf), respectively, are
\begin{align}\label{eq:SR-StDist-def}
H(x)
&\coloneqq
\lim_{t\to+\infty}\Pr(R_{t}^{r}\le x)
\;\;
\text{and}
\;\;
h(x)
\coloneqq
\dfrac{d}{dx}H(x),
\end{align}
where $r\in[0,+\infty)$ is fixed. This probability measure has already been found, e.g., by~\cite{Shiryaev:SMD61,Shiryaev:TPA63}, by~\cite{Pollak+Siegmund:B85}, and more recently also by~\cite{Feinberg+Shiryaev:SD2006,Burnaev+etal:TPA2009,Polunchenko+Sokolov:MCAP2016}, to be the momentless (no moments of orders one and higher) distribution
\begin{align}\label{eq:SR-StDist-answer}
H(x)
&=
e^{-\tfrac{2}{x}}\indicator{x\ge0}
\;\;
\text{and}
\;\;
h(x)
=
\dfrac{2}{x^2}\,e^{-\tfrac{2}{x}}\indicator{x\ge0}
=
\dfrac{2}{x^2}H(x),
\end{align}
which is an extreme-value Fr\'{e}chet-type distribution, and a particular case of the inverse (reciprocal) gamma distribution. See also, e.g.,~\cite{Linetsky:OR2004} and~\cite{Avram+etal:MPRF2013}. As an aside, note that, in view of~\eqref{eq:SR-StDist-answer}, the stationary distribution of the reciprocal of $(R_{t}^{r})_{t\ge0}$ is exponential with mean $1/2$.

However, if all states from a fixed $A>0$ and up inside the process' ``natural'' state space $[0,+\infty)$ are made into absorbing states, then $(R_{t}^{r})_{t\ge0}$ also has a nontrivial probabilistic behavior in the limit as $t\to+\infty$. This behavior is known as quasi-stationarity, and it is characterized by the invariant probability measure whose cdf and pdf, respectively, are
\begin{align}
  Q_A(x) & \coloneqq
  \lim_{t\to+\infty}\Pr(R_{t}^{r}\le x\mid R_{s}^{r}\in[0,A)\;\text{for all}\;0\le s\le t) \label{eq:QSD-def} \\
  \text{and } \; q_A(x) & \coloneqq \dfrac{d}{dx} Q_A(x), \nonumber
\end{align}
where $r\in[0,A)$ is fixed. The existence of this probability measure was formally established, e.g., by~\cite{Pollak+Siegmund:B85}, although one can also infer the same result, e.g., from the earlier seminal work of~\cite{Mandl:CMJ1961}. Moreover, analytic closed-form formulae for both $Q_{A}(x)$ and $q_{A}(x)$ were recently obtained by~\cite{Polunchenko:SA2017a}, apparently for the first time in the literature; see formulae~\eqref{eq:QSD-pdf-answer} and~\eqref{eq:QSD-cdf-answer-W0} in Section~\ref{sec:main-results} below. Recently these formulae were used by~\cite{Polunchenko+Pepelyshev:SP2018} to compute analytically the quasi-stationary distribution's Laplace transform, and then also by~\cite{Li+etal:CommStat2019} to find the quasi-stationary distribution's fractional moment of any real order.
\begin{remark}
The phenomenon of quasi-stationarity is also exhibited by $(R_{t}^{r})_{t\ge0}$ in another case, viz. when all states from 0 up through a fixed $A>0$ inclusive inside the process' ``natural'' state space $[0,+\infty)$ are made into absorbing states, so that the state space of $(R_{t}^{r})_{t\ge0}$ becomes the set $[A,+\infty)$ with absorbtion at the lower end. This case was first investigated in~\cite[Section~7.8.2]{Collet+etal:Book2013}. It was also recently analyzed by~\cite{Polunchenko+etal:TPA2018} who obtained analytically exact closed-form formulae for the corresponding quasi-stationary cdf and pdf.
\end{remark}

The quasi-stationary distribution~\eqref{eq:QSD-def} and the stationary distribution~\eqref{eq:SR-StDist-def} are obviously related: as one would expect, the former converges to the latter as $A\to+\infty$. This was formally shown by~\cite{Pollak+Siegmund:JAP1996}, and not only for the GSR process at hand, but for an entire class of stochastically monotone processes. More specifically, it can be deduced from~\cite{Pollak+Siegmund:JAP1996} that $Q_{A}(x)\ge H(x)$ for any fixed $A>0$ and any $x\ge 0$, and that $\lim_{A\to+\infty} Q_{A}(x)=H(x)$ for any fixed $x\ge0$. The question as to the rate of convergence of $Q_{A}(x)$ down to $H(x)$ as $A\to+\infty$ was recently investigated by~\cite{Li+Polunchenko:SA2020} who showed that $(0<)\;\sup_{x\ge0}[Q_{A}(x)-H(x)]=O(\log(A)/A)$, as $A\to+\infty$; see~\eqref{eq:Q-to-H-unif-conv}. The latter result was arrived at by first obtaining new lower- and upper-bounds for $Q_{A}(x)$, of varying tightness and complexity. This work's main contribution is new lower- and upper-bounds for the pdf $q_{A}(x)$; the bounds for $q_{A}(x)$ can be integrated, and thereby be converted into new bounds for $Q_{A}(x)$. In particular, our new bounds for $q_{A}(x)$ are tight enough to show that $q_{A}(x)$ converges to $h(x)$ as $A\to+\infty$, pointwise, for each fixed $x\ge0$. All of the bounds are obtained explicitly with the aid of the formula for $Q_{A}(x)$ and that for $q_{A}(x)$ latterly offered by~\cite{Polunchenko:SA2017a}, and certain recently discovered monotonicity properties of the modified Bessel $K$ function (of the second kind).

The obtained bounds for $Q_{A}(x)$ and for $q_{A}(x)$ are of importance in quickest change-point detection. Specifically, the process $(R_{t}^{r})_{t\ge0}$ governed by equation~\eqref{eq:Rt_r-def} arises in quickest change-point detection when the aim is to monitor the mean of the process $X_{t}\coloneqq (t-\nu)\indicator{t>\nu}+B_{t}$, where $\nu,t\ge0$, observed ``live''. Since $\EV[X_{t}]=(t-\nu)\indicator{t>\nu}$, it is anticipated that the drift of $(X_{t})_{t\ge0}$ will change from none (zero) to one (per time unit) at time instance $\nu\in[0,+\infty]$ referred to as the change-point. The challenge is that $\nu$ is {\em not} known in advance; in particular $\nu=\infty$ is a possibility, i.e., the drift of $(X_{t})_{t\ge0}$ may remain zero indefinitely and never change. The mean of $(X_{t})_{t\ge0}$ is controlled online by sounding an alarm should (and as soon as) the behavior of $(X_{t})_{t\ge0}$ suggest that possibly $\EV[X_{t}]= t-\nu>0$, i.e., $t>\nu$; if it is not the case, then the alarm is a false one. More concretely, the so-called GSR quickest change-point detection procedure, set up to control the drift of $(X_t)_{t\ge0}$, sounds a false alarm at
\begin{align}\label{eq:T-GSR-def}
\mathcal{S}_{A}^{r}
&\coloneqq
\inf\big\{t\ge0\colon R_{t}^{r}=A\big\}
\;
\text{with}
\;
r\in[0,A)
\;
\text{fixed}
,
\;
\text{and}
\;
\inf\{\varnothing\}=+\infty,
\end{align}
where the constant $A>0$ is selected in advance in accordance with the desired false alarm risk level. Hence $(R_{t}^{r})_{t\ge0}$ is the GSR procedure's detection statistic in the pre-change regime, i.e., for $t\in[0,\nu]$. The definition~\eqref{eq:QSD-def} of the quasi-stationary cdf can now be rewritten as $Q_{A}(x)=\lim_{t\to+\infty}\Pr(R_{t}^{r}\le x|\mathcal{S}_{A}^{r}>t)$.

The GSR procedure, identified in the pre-change regime with the stopping time~\eqref{eq:T-GSR-def}, was proposed by~\cite{Moustakides+etal:SS11} as a headstarted (i.e., more general) version of the classical quasi-Bayesian Shiryaev--Roberts (SR) procedure that emerged from the independent work of Shiryaev~\citeyearpar{Shiryaev:SMD61,Shiryaev:TPA63} and that of Roberts~\citeyearpar{Roberts:T66}. The interest in the GSR procedure (and its variations) is due to its strong (near-) optimality properties. See, e.g.,~\cite{Burnaev:ARSAIM2009},~\cite{Feinberg+Shiryaev:SD2006},~\cite{Burnaev+etal:TPA2009},~\cite{Polunchenko+Tartakovsky:AS10},~\cite{Tartakovsky+Polunchenko:IWAP10},~\cite{Vexler+Gurevich:MESA2011}, and~\cite{Tartakovsky+etal:TPA2012}. For example, it is known that if the GSR procedure's headstart is sampled from the quasi-stationary distribution~\eqref{eq:QSD-def}, then such a randomization of the GSR procedure makes the latter nearly (to within a vanishingly small additive term) minimax in the sense of~\cite{Pollak:AS85}. The idea of such a randomization of the GSR procedure and a proof that the randomized GSR procedure is nearly minimax are due to~\cite{Pollak:AS85} who was concerned with the discrete-time formulation of the problem. For the problem's continuous-time formulation, the same result was established by~\cite{Polunchenko:TPA2017} who heavily relied on the exact closed-form formulae for $Q_{A}(x)$ and $q_{A}(x)$ obtained by~\cite{Polunchenko:SA2017a}, as well as on the quasi-stationary distribution's first two moments, also computed by~\cite{Polunchenko:SA2017a}. The stopping time associated with the randomized GSR procedure is
\begin{align}
  \mathcal{S}_A^Q & \coloneqq \inf\big\{t\ge0\colon R_{t}^{Q}=A\big\} \label{eq:T-SRP-def} \\
  & \text{ with} \; (R_t^Q)_{t\ge0} \; \text{as in~\eqref{eq:Rt_r-def} except} \;
  R_0^Q \propto Q_A(x),  \; \text{and} \; \inf\{\varnothing\} = +\infty , \nonumber
\end{align}
i.e., the initial (at $t=0$) value of the ``original'' $R_{t}^r$ is not a fixed number $r\ge0$, but rather is a random number sampled from the quasi-stationary distribution~\eqref{eq:QSD-def}.

The rest of the paper is four sections. The first one, Section~\ref{sec:nomenclature}, introduces our notation and provides the necessary preliminary background on the special functions needed for our bounds. Section~\ref{sec:preliminary-background} offers a summary of the relevant prior work. The next section, Section~\ref{sec:main-results}, is the paper's main section: this is where we derive our bounds. In Section~\ref{sec:discussion} we illustrate a few applications of our bounds, particularly in quickest change-point detection. Lastly, in Section~\ref{sec:conclusion} we make a few concluding remarks and wrap up the entire paper.

\section{Notation and nomenclature}
\label{sec:nomenclature}

We plan to use the standard mathematical notation. By ``standard'' we mean, e.g., such common nomenclature as $\mathbb{R}$, $\mathbb{C}$, $\mathbb{N}$, $\mathbb{Z}$, the imaginary unit $\mathrm{i}$ defined as the (positive) imaginary ``solution'' of the equation $\mathrm{i}^2=-1$, and so on. More importantly, we will also use the standard notation for a handful of special functions that are to appear repeatedly throughout the sequel. These functions, in their most common notation, are:
\begin{enumerate}
    \setlength{\itemsep}{10pt}
    \setlength{\parskip}{0pt}
    \setlength{\parsep}{0pt}
    \item The Gamma function $\Gamma(z)$, where $z\in\mathbb{C}$, sometimes also regarded as the extension of the factorial to complex numbers, due to the property $\Gamma(n)=(n-1)!$ exhibited for $n\in\mathbb{N}$. See, e.g.,~\cite[Chapter~1]{Bateman+Erdelyi:Book1953v1}.
    \item The (upper-) incomplete Gamma function $\Gamma(a,z)$, where $z\in\mathbb{C}$, defined as
    \begin{align*} 
    \Gamma(a,z)
    &\coloneqq
    \int_{z}^{+\infty} y^{a}\, e^{-y}\dfrac{dy}{y},
    \end{align*}
    with no restriction on the integration path. See~\cite[Chapter~9]{Bateman+Erdelyi:Book1953v2}. The ``complete'' Gamma function $\Gamma(z)$ introduced earlier is a special case of the incomplete Gamma function $\Gamma(a,z)$ because $\Gamma(0,z)=\Gamma(z)$.
    \item The exponential integral function $\Ei(x)$, where $x\in\mathbb{R}\backslash\{0\}$, defined as
    \begin{align}\label{eq:Ei-func-def}
    \Ei(x)
    &\coloneqq%
    \begin{cases}
    -\displaystyle\int_{-x}^{+\infty} e^{-y}\dfrac{dy}{y},&\text{if $x<0$;}\\[4mm]
    -\lim_{\varepsilon\to+0}\left[\displaystyle\int_{-x}^{-\varepsilon}e^{-y}\dfrac{dy}{y}+\displaystyle\int_{\varepsilon}^{+\infty}e^{-y}\dfrac{dy}{y}\right],&\text{if $x>0$},\\[2mm]
    \end{cases}
    \end{align}
    with a singularity at $x=0$. Its basic properties are summarized, e.g., in~\cite[Chapter~5]{Abramowitz+Stegun:Handbook1964}. More specifically, we will need the function $\E1(x)\coloneqq-\Ei(-x)$ with $x>0$, i.e.,
    \begin{align}\label{eq:E1-func-def}
    \E1(x)
    &\coloneqq\Gamma(0,x)=
    \int_{x}^{+\infty} e^{-y}\,\dfrac{dy}{y},\; x>0;
    \end{align}
    see also, e.g.,~\cite[Chapter~5]{Abramowitz+Stegun:Handbook1964}.
    \item The Whittaker $M$ and $W$ functions, traditionally denoted, respectively, as $M_{a,b}(z)$ and $W_{a,b}(z)$, where $a,b,z\in\mathbb{C}$. These functions were introduced by Whittaker~\citeyearpar{Whittaker:BAMS1904} as the fundamental solutions to the Whittaker differential equation. See, e.g.,~\cite{Slater:Book1960} and~\cite{Buchholz:Book1969}.
    \item The modified Bessel functions of the first and second kinds, conventionally denoted, respectively, as $I_{a}(z)$ and $K_{a}(z)$, where $a,z\in\mathbb{C}$; the index $a$ is referred to as the function's order. See~\cite[Chapter~7]{Bateman+Erdelyi:Book1953v2}. These functions form a set of fundamental solutions to the modified Bessel differential equation. The modified Bessel $K$ function is also known as the MacDonald function.
\end{enumerate}


\section{Preliminary background on the quasi-stationary distribution}
\label{sec:preliminary-background}

The quasi-stationary distribution's pdf and cdf defined by~\eqref{eq:QSD-def} can both be expressed analytically and in closed form; see~\cite{Polunchenko:SA2017a}. We now recall the expressions, as they will be key to establishing our main results in the next section. Specifically, it can be deduced from~\cite[Theorem~3.1]{Polunchenko:SA2017a} that if $A>0$ is fixed and $\lambda\equiv\lambda_A>0$ is the smallest (positive) solution of the equation
\begin{align}\label{eq:lambda-eqn}
W_{1,\tfrac{1}{2}\xi(\lambda)}\left(\dfrac{2}{A}\right)
&=
0,
\end{align}
where
\begin{align}\label{eq:xi-def}
\xi(\lambda)
&\coloneqq
\sqrt{1-8\lambda}
\;\;
\text{so that}
\;\;
\lambda
=
\dfrac{1}{8}\left(1-\big[\xi(\lambda)\big]^2\right),
\end{align}
then the quasi-stationary pdf is given by
\begin{align}\label{eq:QSD-pdf-answer}
q_{A}(x)
&=
\dfrac{e^{-\tfrac{1}{x}}\,\dfrac{1}{x}\,W_{1,\tfrac{1}{2}\xi(\lambda)}\left(\dfrac{2}{x}\right)}{e^{-\tfrac{1}{A}}\,W_{0,\tfrac{1}{2}\xi(\lambda)}\left(\dfrac{2}{A}\right)}\indicator{x\in[0,A]},
\end{align}
and the respective cdf is given either by
\begin{align}\label{eq:QSD-cdf-answer-W0}
Q_{A}(x)
&=
\begin{cases}
1,&\;\text{if $x\ge A$;}\\[2mm]
\dfrac{e^{-\tfrac{1}{x}}\,W_{0,\tfrac{1}{2}\xi(\lambda)}\left(\dfrac{2}{x}\right)}{e^{-\tfrac{1}{A}}\,W_{0,\tfrac{1}{2}\xi(\lambda)}\left(\dfrac{2}{A}\right)},&\;\text{if $x\in[0,A)$;}\\[8mm]
0,&\;\text{otherwise},
\end{cases}
\end{align}
or, equivalently, by
\begin{align}\label{eq:QSD-cdf-answer-K}
Q_{A}(x)
&=
\begin{cases}
1,&\;\text{if $x\ge A$;}\\[2mm]
\sqrt{\dfrac{A}{x}}\,\dfrac{e^{-\tfrac{1}{x}}\,K_{\tfrac{1}{2}\xi(\lambda)}\left(\dfrac{1}{x}\right)}{e^{-\tfrac{1}{A}}\,K_{\tfrac{1}{2}\xi(\lambda)}\left(\dfrac{1}{A}\right)},&\;\text{if $x\in[0,A)$;}\\[8mm]
0,&\;\text{otherwise},
\end{cases}
\end{align}
because
\begin{align}\label{eq:Whit0-BesselK-id}
W_{0,b}(z)
&=
\sqrt{\dfrac{z}{\pi}}\, K_{b}\left(\dfrac{z}{2}\right),
\end{align}
which is~\cite[Identity~9.6.48,~p.~377]{Abramowitz+Stegun:Handbook1964}. Formula~\eqref{eq:QSD-cdf-answer-W0} is a special case of~\cite[Formula~(3.11),~p.~134]{Polunchenko:SA2017a}. Formula~\eqref{eq:QSD-cdf-answer-K}, in turn, is precisely~\cite[Formula~(3.8),~p.~220]{Li+Polunchenko:SA2020}. Observe also that~\eqref{eq:lambda-eqn}, \eqref{eq:xi-def}, and~\eqref{eq:QSD-pdf-answer} together yield $q_{A}(A)=0$ for any $A>0$.

The pdf $q_{A}(x)$ is a ``singularity-free'', bounded function in both $x\in\mathbb{R}$ as well as $A>0$, even in the limit as $A\to+\infty$; so is the cdf $Q_{A}(x)$, of course. This is due to certain analytic properties of the Whittaker $W$ function on the right of~\eqref{eq:QSD-pdf-answer}. As a result, one can, for example, differentiate and/or integrate $q_{A}(x)$ in any order without any issues. Another relevant consequence is~\cite[Lemma~3.1]{Polunchenko:SA2017a} whereby
\begin{align*}
\lim_{x\to0+}\left[q_{A}(x)\right]=0,
\;
\text{for any {\em fixed}}
\;
A>0,
\end{align*}
and
\begin{align*}
\lim_{x\to0+}\left\{\dfrac{\partial}{\partial x}\left[q_{A}(x)\right]\right\}=0,
\;
\text{for any {\em fixed}}
\;
A>0,
\end{align*}
and both limits will be used in the sequel.

Formulae~\eqref{eq:QSD-pdf-answer} and~\eqref{eq:QSD-cdf-answer-W0} stem from the solution of a certain Sturm--Liouville problem, and $\lambda$ is the smallest positive eigenvalue of the corresponding Sturm--Liouville operator; if the Sturm--Liouville operator is negated, as was done by~\cite{Polunchenko:SA2017a}, then $\lambda$ becomes the operator's largest {\em negative} eigenvalue.
\begin{remark}\label{rem:xi-symmetry}
The definition~\eqref{eq:xi-def} of $\xi(\lambda)$ can actually be changed to $\xi(\lambda)\coloneqq -\sqrt{1-8\lambda}$ with no effect whatsoever on either equation~\eqref{eq:lambda-eqn}, or formulae~\eqref{eq:QSD-pdf-answer} and~\eqref{eq:QSD-cdf-answer-W0}, i.e., all three are invariant with respect to the sign of $\xi(\lambda)$. This was previously pointed out by~\cite{Polunchenko:SA2017a}, and the reason for this $\xi(\lambda)$-symmetry is because equation~\eqref{eq:lambda-eqn} and formulae~\eqref{eq:QSD-pdf-answer} and~\eqref{eq:QSD-cdf-answer-W0} each have $\xi(\lambda)$ present only as (double) the second index of the corresponding Whittaker $W$ function or functions involved, and the Whittaker $W$ function in general is known (see, e.g.,~\cite[Identity~(19),~p.~19]{Buchholz:Book1969}) to be an even function of its second index, i.e., $W_{a,b}(z)=W_{a,-b}(z)$.
\end{remark}

It is evident that equation~\eqref{eq:lambda-eqn} is a key component of formulae~\eqref{eq:QSD-pdf-answer} and~\eqref{eq:QSD-cdf-answer-W0}, and consequently, of all of the characteristics of the quasi-stationary distribution as well. As a transcendental equation, it can only be solved numerically, although to within any desired accuracy; see, e.g.~\cite{Linetsky:OR2004,Polunchenko:SA2016,Polunchenko:SA2017a,Polunchenko:SA2017b}. Yet, it is known (see, e.g.,~\citealt{Linetsky:OR2004} and~\citealt{Polunchenko:SA2016}) that for any fixed $A>0$, the equation has countably many simple solutions $0<\lambda_1<\lambda_2<\lambda_3<\cdots$, such that $\lim_{k\to+\infty}\lambda_k=+\infty$. All of them, of course, do depend on $A$, but since we are interested only in the smallest one, we shall use either the ``short'' notation $\lambda$, or the more explicit $\lambda_A$ to emphasize the dependence on $A$. It was shown by~\cite{Polunchenko:SA2017a} that $\lambda_A$ is a monotonically decreasing function of $A$, and such that
\begin{align}\label{eq:lambda-dbl-ineq}
\dfrac{1}{A}
+
\dfrac{1}{A(1+A)}
&<
\lambda_{A}
<
\dfrac{1}{A}
+
\dfrac{1+\sqrt{4A+1}}{2A^2}
,
\;\;
\text{for any}
\;\;
A>0,
\end{align}
whence $\lim_{A\to+\infty}\lambda_A=0$, and more specifically $\lambda_A=A^{-1}+O(A^{-3/2})$; cf.~\cite[p.~136 and Lemma~3.3]{Polunchenko:SA2017a}. See also~\cite{Polunchenko+Pepelyshev:SP2018} for a discussion of potential ways to improve the foregoing double inequality.
\begin{remark}\label{rem:xi-complex-real}
Since $\lambda\equiv\lambda_{A}$ is monotonically decreasing in $A$, and such that $\lim_{A\to+\infty}\lambda_{A}=0$, one can conclude from~\eqref{eq:xi-def} that $\xi(\lambda_{A})$, for any finite $A>0$, is either \begin{inparaenum}[\itshape(a)]\item purely imaginary (i.e., $\xi(\lambda)=\mathrm{i}\alpha$ where $\mathrm{i}\coloneqq\sqrt{-1}$ and $\alpha\in\mathbb{R}$) if $A$ is sufficiently small, or \item purely real and between 0 inclusive and 1 exclusive (i.e., $0\le \xi(\lambda)<1$) otherwise\end{inparaenum}. The borderline case is when $\xi(\lambda)=0$, i.e., when $\lambda_{A}=1/8$, and the corresponding critical value of $A$ is the solution $\tilde{A}>0$ of the equation
\begin{align*} 
  W_{1,0}\big( 2 / \tilde{A} \big)
  & = 0, \;\; \text{so that} \;\; \tilde{A}\approx10.240465,
\end{align*}
as can be established by a basic numerical calculation. Hence, if $A<\tilde{A}\approx10.240465$, then $\lambda_{A}>1/8$ so that $\xi(\lambda)$ is purely imaginary; otherwise, if $A\ge\tilde{A}\approx10.240465$, then $\lambda_{A}\in(0,1/8]$ so that $\xi(\lambda)$ is purely real and such that $\xi(\lambda)\in[0,1)$ with $\lim_{A\to+\infty}\xi(\lambda_{A})=1$.
\end{remark}

To numerically evaluate $q_A(x)$ and/or $Q_A(x)$ one can use such software packages as R, Python, or {\it Wolfram Mathematica}. In the sequel we will use the latter (specifically, Mathematica's routine called \texttt{BesselK}).

It was recently shown by~\cite{Li+Polunchenko:SA2020} that
\begin{align}
  0 \le \sup_{x\in\mathbb{R}}\left[Q_{A}(x)-H(x)\right] & = \sup_{x\in\mathbb{R}}\abs{Q_{A}(x)-H(x)}
  = O\left(\dfrac{\log(A)}{A}\right), \label{eq:Q-to-H-unif-conv} \\
  & \quad \text{as} \; A\to+\infty, \nonumber
\end{align}
i.e., $Q_A(x)$ converges to $H(x)$ uniformly in $x$, as $A\to+\infty$. The proof makes use of the double-inequality~\eqref{eq:lambda-dbl-ineq} and formula~\eqref{eq:QSD-cdf-answer-K}, along with certain latest monotonicity properties of the Bessel $K$ function and functionals thereof.


\section{Bounds for the quasi-stationary distribution}
\label{sec:main-results}

We are now in a position to start deriving our bounds for the quasi-stationary distribution. To that end, formula~\eqref{eq:QSD-pdf-answer} for the pdf $q_{A}(x)$ is poorly suited for our purposes: the Whittaker $W$ function involved in~\eqref{eq:QSD-pdf-answer} is a somewhat exotic special function (with $1$ as the first index, anyway), and has not received much attention in the literature. Nevertheless, it {\em is} possible to rid~\eqref{eq:QSD-pdf-answer} of the Whittaker $W$ function, and instead express $q_{A}(x)$ entirely in terms of the modified Bessel $K$ function---a much more extensively studied special function. Better yet, alternative, well-suited for our purposes expressions for $q_{A}(x)$ are many. For example, the pdf $q_{A}(x)$ is known to satisfy two second-order ordinary differential equations: one called the Kolmogorov forward and the other called the Kolmogorov backward equation. The equations are related, and either one alone, subject to the appropriate boundary and regularity conditions, uniquely defines $q_{A}(x)$. The boundary condition is $q_{A}(A)=0$, which is a Dirichlet-type condition effectively saying that the GSR process is ``killed'' at $x=A$. The regularity condition is that $q_{A}(x)$ as a pdf must integrate to unity over $[0,A]$. Formula~\eqref{eq:QSD-pdf-answer} for $q_{A}(x)$ was arrived at by~\cite{Polunchenko:SA2017a} by solving explicitly the forward equation, which is
\begin{align*}
\dfrac{1}{2}\dfrac{\partial^2}{\partial y^2}\left[y^2 q_{A}(y)\right]-\dfrac{\partial}{\partial y}\left[q_{A}(y)\right]
&=
-\lambda q_{A}(y),
\;
y\in[0,A],
\;
A>0,
\end{align*}
where $\lambda\equiv\lambda_{A}$ is determined by~\eqref{eq:lambda-eqn}. If we now integrate this equation with respect to $y$ from 0 up through $x\in[0,A]$, then we will get
\begin{align}\label{eq:QSD-pdf-qA-via-cdf-QA-v1}
q_{A}(x)
&=
\dfrac{2}{x^2}\left\{Q_{A}(x)-\lambda\int_{0}^{x}Q_{A}(t)\,dt\right\}\indicator{x\in[0,A]},
\;
x\in\mathbb{R},
\;
A>0,
\end{align}
or equivalently
\begin{align}\label{eq:QSD-pdf-qA-via-cdf-QA-v2}
q_{A}(x)
&=
\dfrac{2}{x^2}\left\{Q_{A}(x)-1+\lambda\int_{x}^{A}Q_{A}(t)\,dt\right\}\indicator{x\in[0,A]},
\;
x\in\mathbb{R},
\;
A>0,
\end{align}
inasmuch as $q_{A}(A)=0$ and $Q_{A}(A)=1$, for any $A>0$, due to~\eqref{eq:lambda-eqn},~\eqref{eq:QSD-pdf-answer}, and~\eqref{eq:QSD-cdf-answer-W0}.

Parenthetically, we note the curious identity
\begin{align}\label{eq:QA-int-lambda-link}
\int_{0}^{A}Q_{A}(t)\,dt
&=
\dfrac{1}{\lambda_{A}},
\;
A>0,
\end{align}
which comes directly from~\eqref{eq:QSD-pdf-qA-via-cdf-QA-v1} with $x=A$, because again $q_{A}(A)=0$ and $Q_{A}(A)=1$ for any $A>0$. Identity~\eqref{eq:QA-int-lambda-link} is noteworthy because, with its help, any lower- or upper-bound for $Q_{A}(x)$ can be converted into an upper- or, respectively, a lower-bound for $\lambda$, and the behavior of $\lambda_{A}$ as a function of $A>0$ is of importance in quickest change-point detection; see, e.g.,~\cite{Polunchenko:TPA2017}. For example, by virtue of~\eqref{eq:QA-int-lambda-link}, one can get the trivial inequality $\lambda>1/A$, $A>0$, from the trivial bound $Q_{A}(x)\le 1$ valid for all $x\in\mathbb{R}$ and any $A>0$. Getting sharper bounds, such as~\eqref{eq:lambda-dbl-ineq} or even better, is also possible, and this question will be explored in the next section.


%

Now, with virtually no effort one can get a pair of lower- and a pair of upper-bounds for $q_{A}(x)$ directly from~\eqref{eq:QSD-pdf-qA-via-cdf-QA-v1} and~\eqref{eq:QSD-pdf-qA-via-cdf-QA-v2}. Specifically, since $Q_{A}(x)$ is strictly increasing in $x$ for any $A>0$, and such that $0\le Q_{A}(x)\le 1$ for all $x\in\mathbb{R}$ for any $A>0$, from~\eqref{eq:QSD-pdf-qA-via-cdf-QA-v1} we get
\begin{align}\label{eq:QSD-pdf-qA-lwrbnd1-def}
l_{A}^{(1)}(x)
&\coloneqq
\dfrac{2}{x^2}Q_{A}(x)[1-\lambda x]\indicator{x\in[0,A]}
\le
q_{A}(x),
\;
x\in\mathbb{R},
\;
A>0,
\end{align}
but
\begin{align}\label{eq:QSD-pdf-qA-uprbnd1-def}
q_{A}(x)
&\le
\dfrac{2}{x^2}Q_{A}(x)\indicator{x\in[0,A]}\eqqcolon u_{A}^{(1)}(x),
\;
x\in\mathbb{R},
\;
A>0,
\end{align}
and, likewise, from~\eqref{eq:QSD-pdf-qA-via-cdf-QA-v2} we get
\begin{align*}  
l_{A}^{(2)}(x)
&\coloneqq
\dfrac{2}{x^2}\left[Q_{A}(x)(1+\lambda A-\lambda x)-1\right]\indicator{x\in[0,A]}
\le
q_{A}(x),
\;
x\in\mathbb{R},
\;
A>0,
\end{align*}
but
\begin{align*} 
q_{A}(x)
&\le
\dfrac{2}{x^2}\left[Q_{A}(x)+\lambda A-1-\lambda x\right]\indicator{x\in[0,A]}\eqqcolon u_{A}^{(2)}(x),
\;
x\in\mathbb{R},
\;
A>0.
\end{align*}

The pair $(l_{A}^{(1)}(x),u_{A}^{(1)}(x))$ and the pair $(l_{A}^{(2)}(x),u_{A}^{(2)}(x))$ complement each other. Specifically, on the one hand, for the first pair we have $\lim_{x\to0+}\left[l_{A}^{(1)}(x)\right]=\lim_{x\to0+}\left[u_{A}^{(1)}(x)\right]=0$, for any fixed $A>0$, but $l_{A}^{(1)}(A)=2(1-\lambda A)/A^2<0$, due to the left half of~\eqref{eq:lambda-dbl-ineq}, and $u_{A}^{(1)}(A)=2/A^2>0$, for any $A>0$. That is, the bounds $l_{A}^{(1)}(x)$ and $u_{A}^{(1)}(x)$ are loose if $x$ is sufficiently close to $A$, although note that $\lim_{A\to+\infty}\left[l_{A}^{(1)}(A)\right]=0$ and $\lim_{A\to+\infty}\left[u_{A}^{(1)}(A)\right]=0$, because of~\eqref{eq:lambda-dbl-ineq}. On the other hand, by contrast, for the second pair of bounds we have $l_{A}^{(2)}(A)=u_{A}^{(2)}(A)=0$ for any $A>0$, but $\lim_{x\to0+}\left[l_{A}^{(2)}(x)\right]=-\infty$ and $\lim_{x\to0+}\left[u_{A}^{(2)}(x)\right]=+\infty$, for any fixed $A>0$. That is, the bounds $l_{A}^{(2)}(x)$ and $u_{A}^{(2)}(x)$ are (too) loose if $x$ is close to $0$, no matter whether $A$ is small or large. This suggests that one can consider
\begin{align*}
\max\left\{l_{A}^{(1)}(x),l_{A}^{(2)}(x)\right\}
&\le
q_{A}(x)
\le
\min\left\{u_{A}^{(1)}(x),u_{A}^{(2)}(x)\right\},
\;
x\in\mathbb{R},
\;
A>0,
\end{align*}
although, for simplicity's sake, we shall just stick to $l_{A}^{(1)}(x)$ and $u_{A}^{(1)}(x)$, which, unlike $l_{A}^{(2)}(x)$ and $u_{A}^{(2)}(x)$, are singularity-free, and get sharp for all $x$ asymptotically, as $A\to+\infty$.

\begin{figure}[h!]
    \centering
    \begin{subfigure}{0.48\textwidth}
        \centering
        \includegraphics[width=\linewidth]{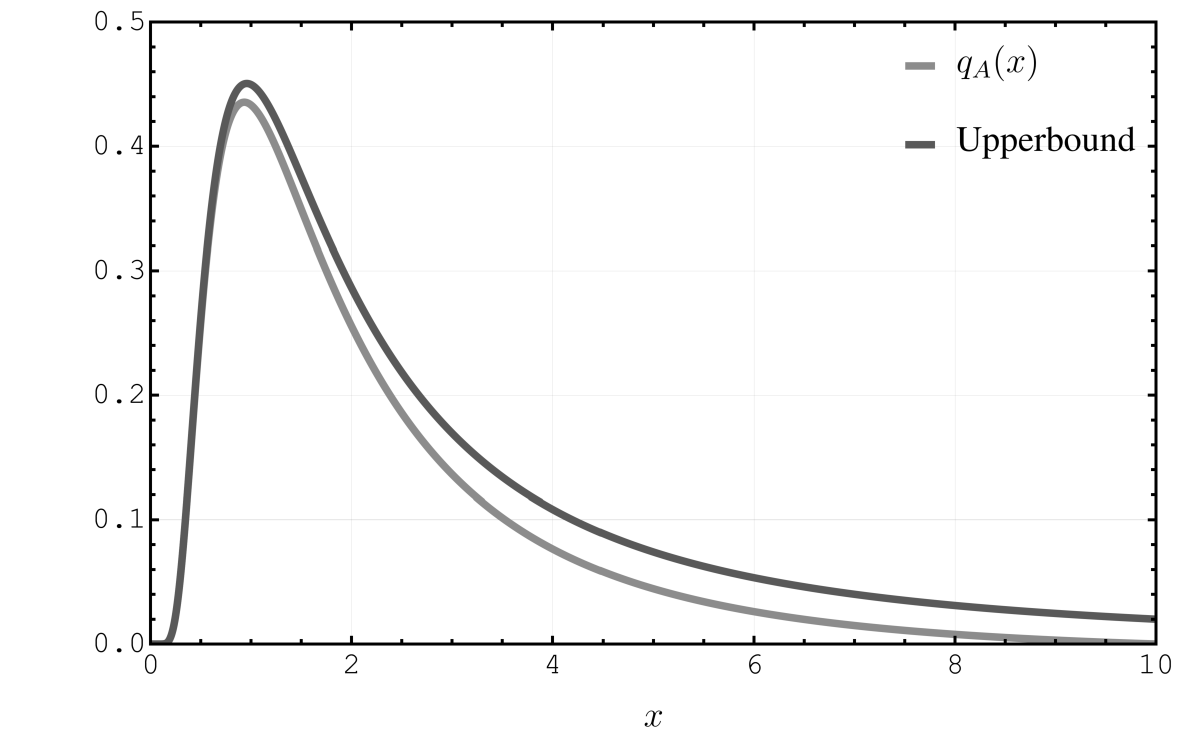}
        \caption{$q_{A}(x)$ and $u_{A}^{(1)}(x)$.}
        \label{fig:QST-pdf-uprbnd1-A10}
    \end{subfigure}
    \hspace*{\fill}
    \begin{subfigure}{0.48\textwidth}
        \centering
        \includegraphics[width=\linewidth]{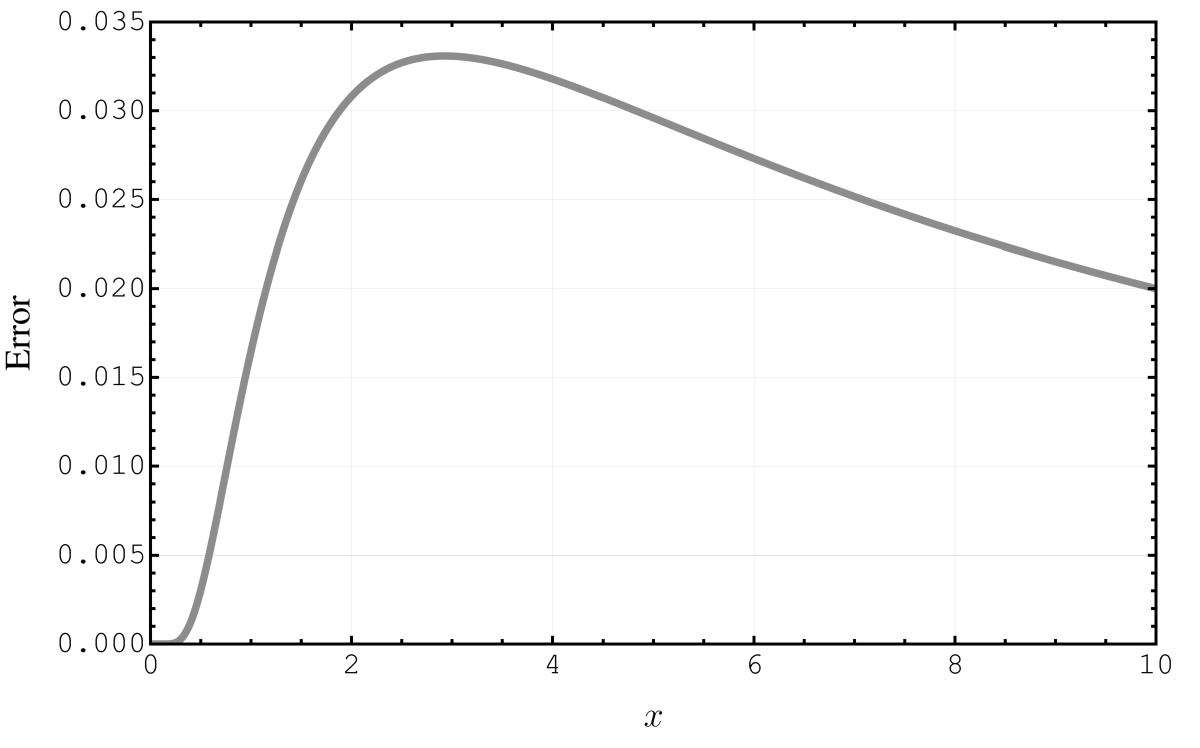}
        \caption{Corresponding upperbound error.}
        \label{fig:QST-pdf-uprbnd1-err-A10}
    \end{subfigure}
    \caption{Quasi-stationary distribution's pdf, $q_{A}(x)$, its upperbound $u_A^{(1)}(x)$, and the corresponding error---all as functions of $x\in[0,A]$ for $A=10$.}
    \label{fig:QST-pdf-uprbnd1-perf-A10}
\end{figure}

\begin{figure}[h!]
    \centering
    \begin{subfigure}{0.48\textwidth}
        \centering
        \includegraphics[width=\linewidth]{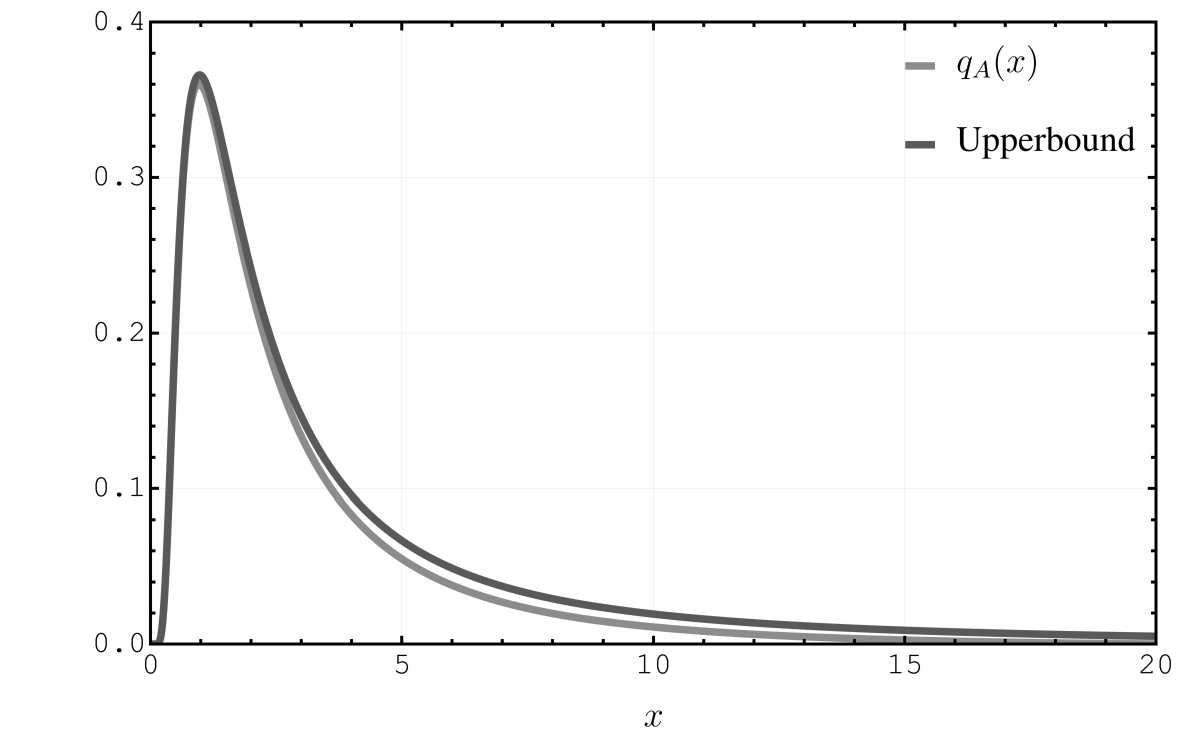}
        \caption{$q_{A}(x)$ and $u_A^{(1)}(x)$.}
        \label{fig:QST-pdf-uprbnd1-A20}
    \end{subfigure}
    \hspace*{\fill}
    \begin{subfigure}{0.48\textwidth}
        \centering
        \includegraphics[width=\linewidth]{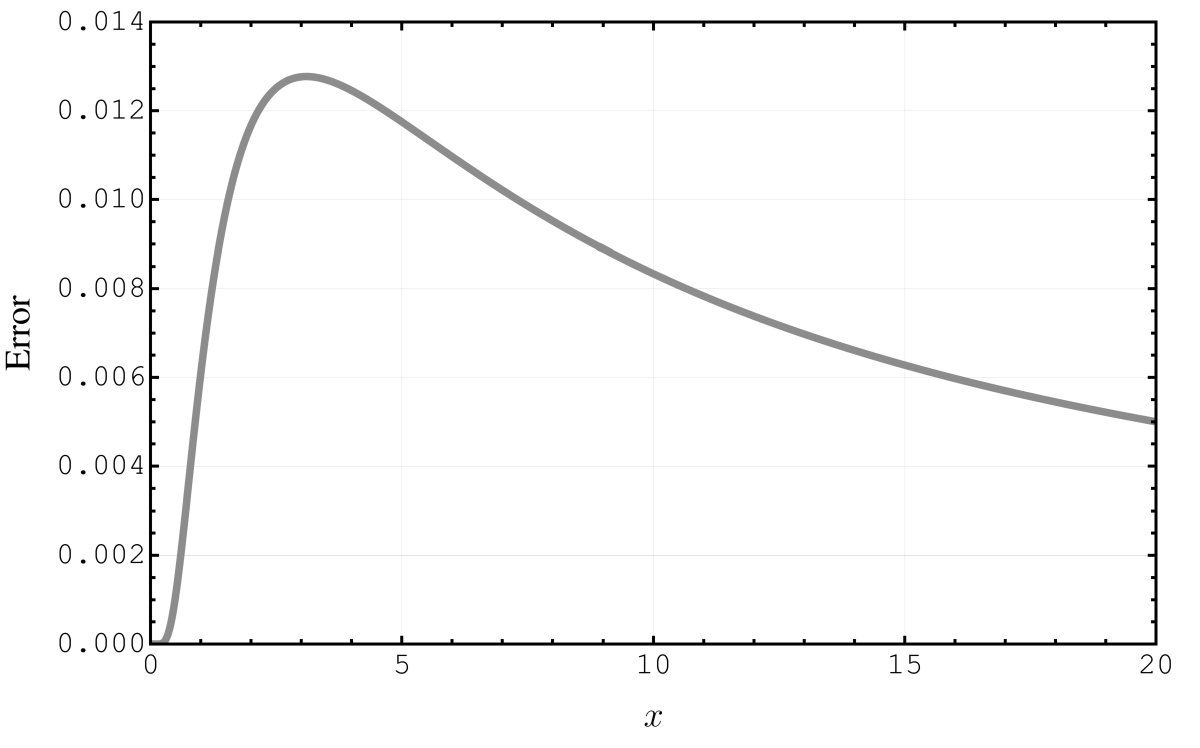}
        \caption{Corresponding upperbound error.}
        \label{fig:QST-pdf-uprbnd1-err-A20}
    \end{subfigure}
    \caption{Quasi-stationary distribution's pdf, $q_{A}(x)$, its upperbound $u_A^{(1)}(x)$, and the corresponding error---all as functions of $x\in[0,A]$ for $A=20$.}
    \label{fig:QST-pdf-uprbnd1-perf-A20}
\end{figure}

\begin{figure}[h!]
    \centering
    \begin{subfigure}{0.48\textwidth}
        \centering
        \includegraphics[width=\linewidth]{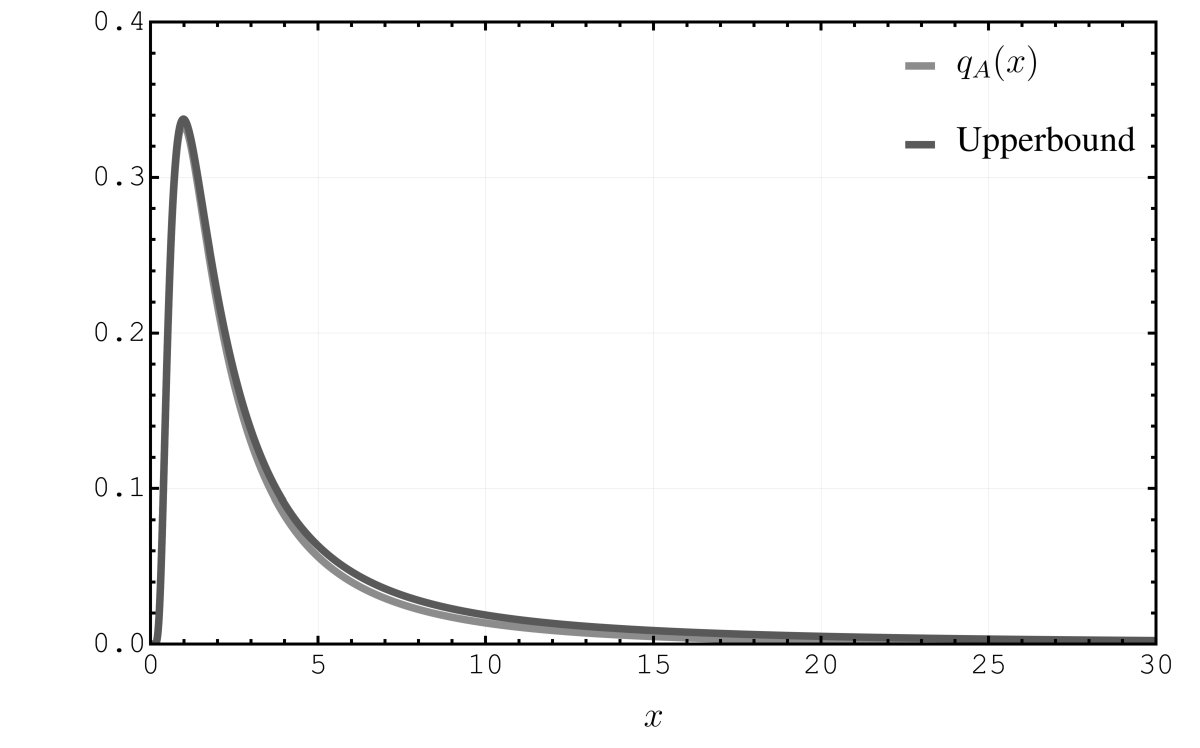}
        \caption{$q_{A}(x)$ and $u_A^{(1)}(x)$.}
        \label{fig:QST-pdf-uprbnd1-A30}
    \end{subfigure}
    \hspace*{\fill}
    \begin{subfigure}{0.48\textwidth}
        \centering
        \includegraphics[width=\linewidth]{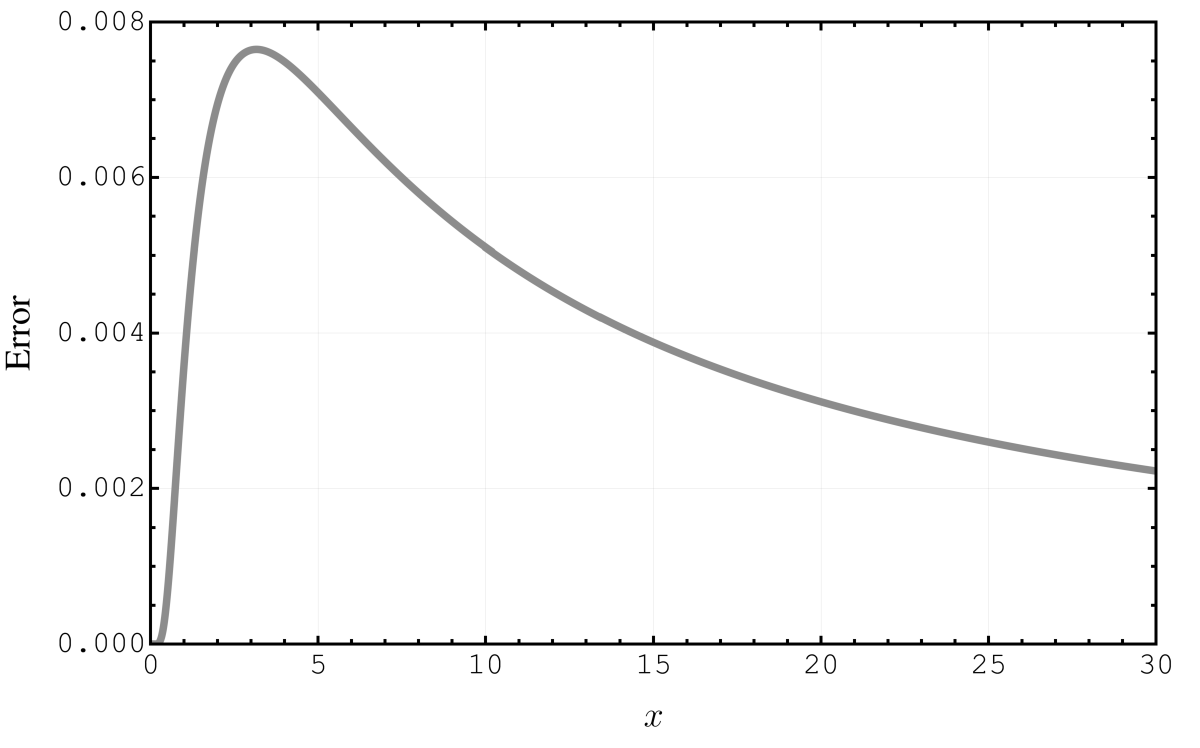}
        \caption{Corresponding upperbound error.}
        \label{fig:QST-pdf-uprbnd1-err-A30}
    \end{subfigure}
    \caption{Quasi-stationary distribution's pdf, $q_{A}(x)$, its upperbound $u_A^{(1)}(x)$, and the corresponding error---all as functions of $x\in[0,A]$ for $A=30$.}
    \label{fig:QST-pdf-uprbnd1-perf-A30}
\end{figure}

The lowerbound $l_{A}^{(1)}(x)$ and the upperbound $u_{A}^{(1)}(x)$---as ``simple'' as they may seem---are actually fairly tight, for all $x\in[0,A]$, even if $A$ is low. This is best demonstrated numerically, with the aid of the {\it Wolfram Mathematica} script prepared by~\cite{Polunchenko:SA2017a}: the script can compute $\lambda$ to within five hundred (!) decimal places of accuracy, so that both $q_{A}(x)$ and $Q_{A}(x)$ can be computed almost precisely. For example, Figures~\ref{fig:QST-pdf-uprbnd1-perf-A10},~\ref{fig:QST-pdf-uprbnd1-perf-A20}, and~\ref{fig:QST-pdf-uprbnd1-perf-A30} show the upperbound $u_A^{(1)}(x)$ and the actual pdf $q_A(x)$ as functions of $x\in[0,A]$ for $A$ set as low as $10$, $20$, and $30$, respectively. It can be seen from the figures that the discrepancy between $u_{A}^{(1)}(x)$ and $q_{A}(x)$ is fairly small, for all $x\in[0,A]$, and rapidly gets even smaller as $A$ increases.

\begin{figure}[h!]
    \centering
    \begin{subfigure}{0.48\textwidth}
        \centering
        \includegraphics[width=\linewidth]{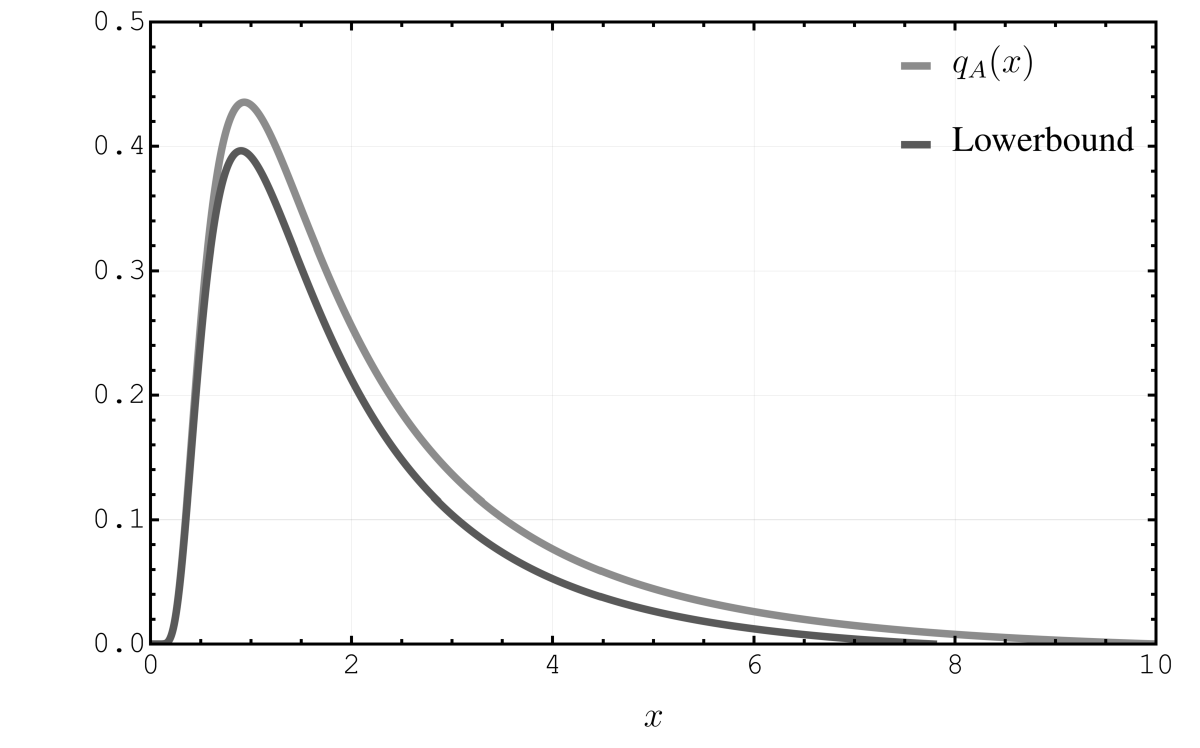}
        \caption{$q_{A}(x)$ and $l_{A}^{(1)}(x)$.}
        \label{fig:QST-pdf-lwrbnd1-A10}
    \end{subfigure}
    \hspace*{\fill}
    \begin{subfigure}{0.48\textwidth}
        \centering
        \includegraphics[width=\linewidth]{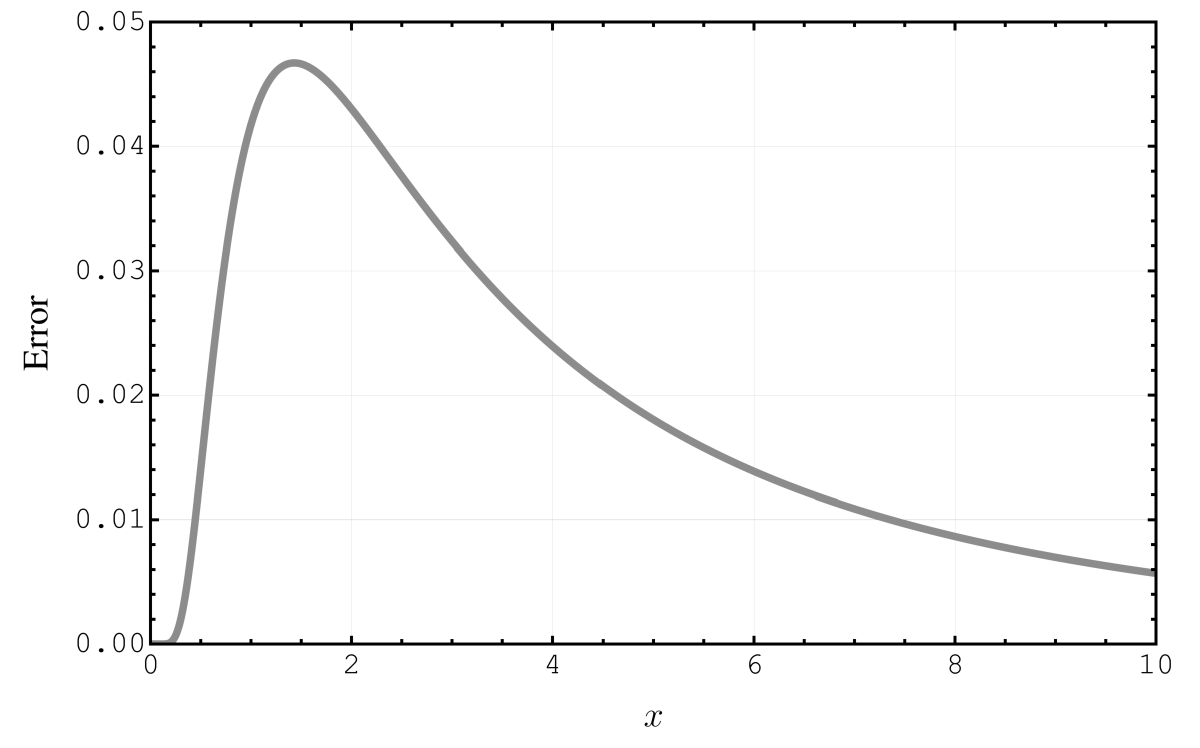}
        \caption{Corresponding lowerbound error.}
        \label{fig:QST-pdf-uprbnd1-err-A10}
    \end{subfigure}
    \caption{Quasi-stationary distribution's pdf, $q_{A}(x)$, its lowerbound $l_A^{(1)}(x)$, and the corresponding error---all as functions of $x\in[0,A]$ for $A=10$.}
    \label{fig:QST-pdf-lwrbnd1-perf-A10}
\end{figure}

\begin{figure}[h!]
    \centering
    \begin{subfigure}{0.48\textwidth}
        \centering
        \includegraphics[width=\linewidth]{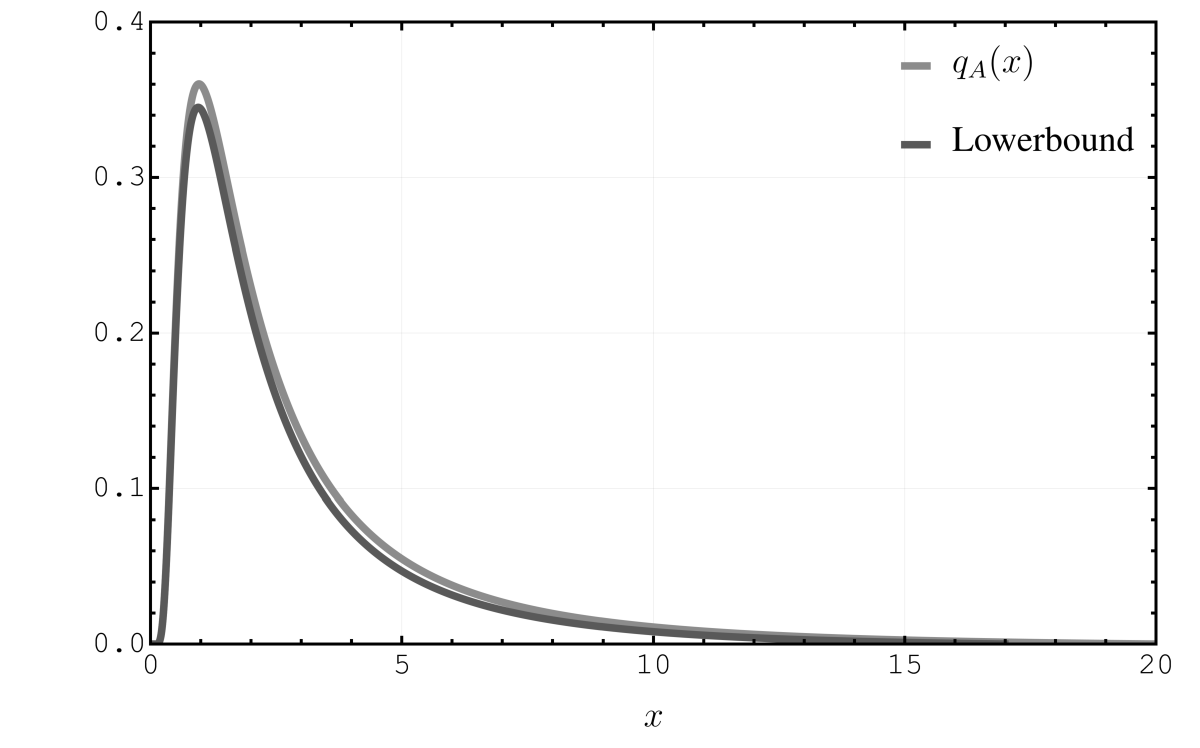}
        \caption{$q_{A}(x)$ and $l_{A}^{(1)}(x)$.}
        \label{fig:QST-pdf-lwrbnd1-A20}
    \end{subfigure}
    \hspace*{\fill}
    \begin{subfigure}{0.48\textwidth}
        \centering
        \includegraphics[width=\linewidth]{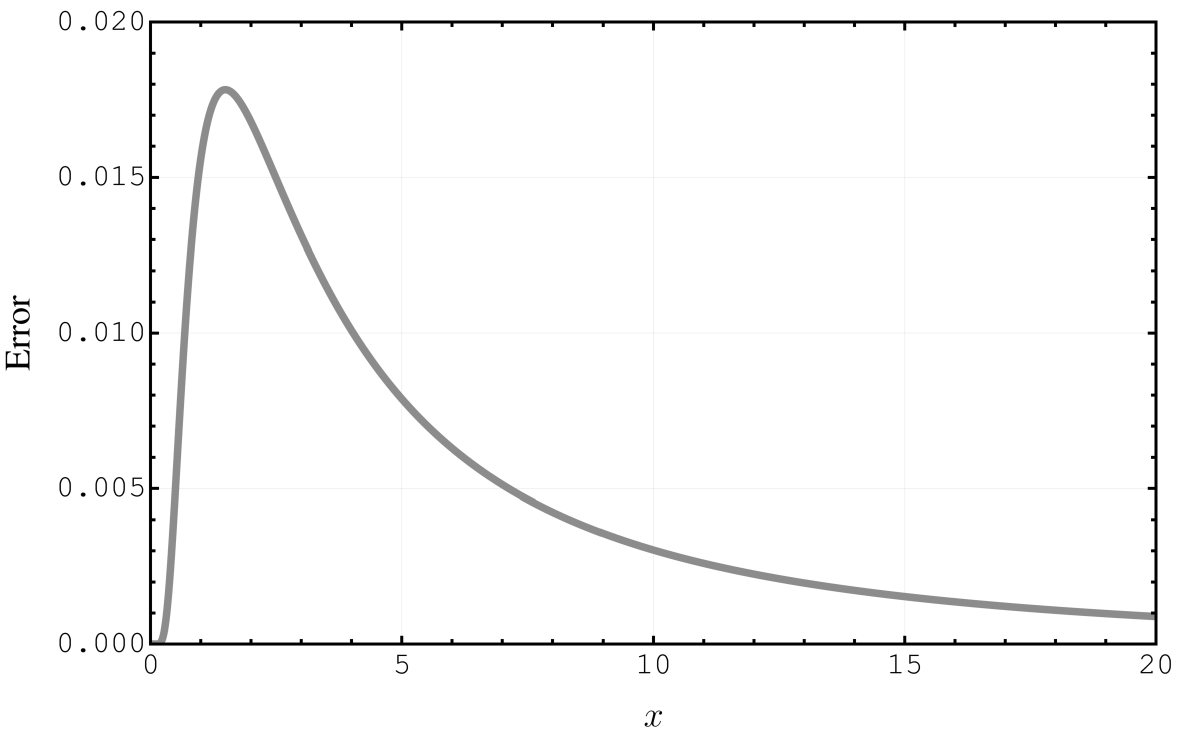}
        \caption{Corresponding lowerbound error.}
        \label{fig:QST-pdf-uprbnd1-err-A20}
    \end{subfigure}
    \caption{Quasi-stationary distribution's pdf, $q_{A}(x)$, its lowerbound $l_A^{(1)}(x)$, and the corresponding error---all as functions of $x\in[0,A]$ for $A=20$.}
    \label{fig:QST-pdf-lwrbnd1-perf-A20}
\end{figure}

\begin{figure}[h!]
    \centering
    \begin{subfigure}{0.48\textwidth}
        \centering
        \includegraphics[width=\linewidth]{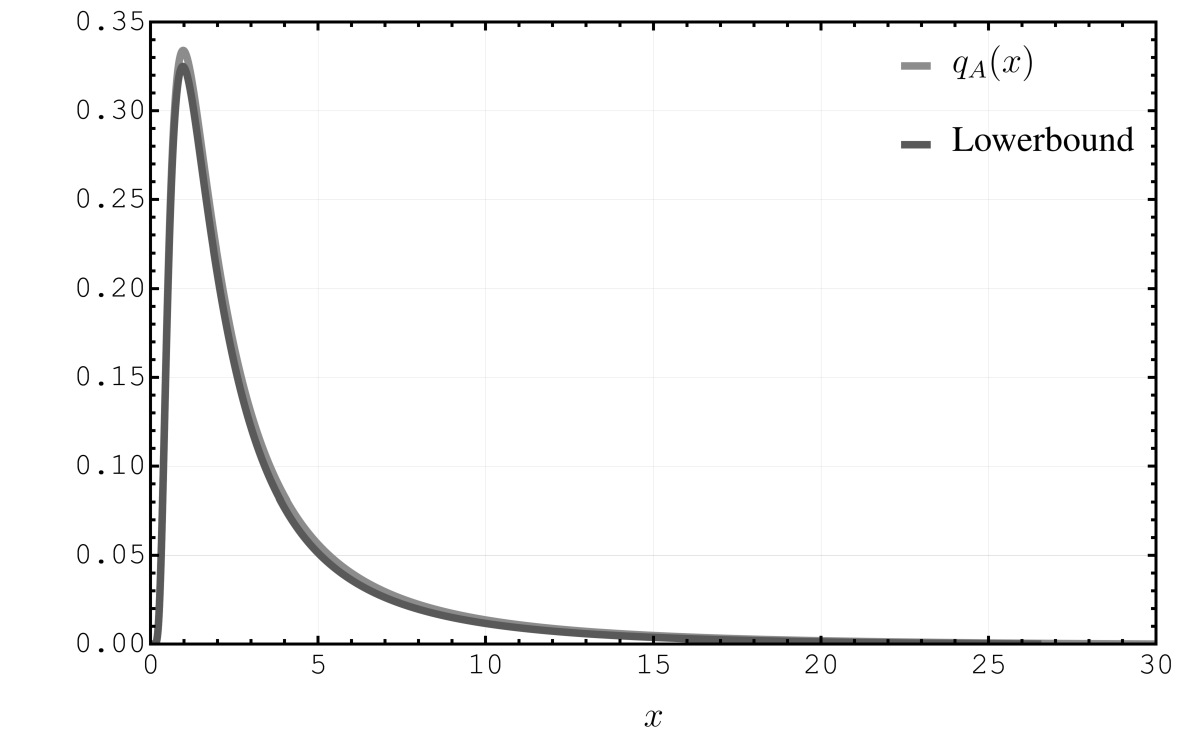}
        \caption{$q_{A}(x)$ and $l_{A}^{(1)}(x)$.}
        \label{fig:QST-pdf-lwrbnd1-A30}
    \end{subfigure}
    \hspace*{\fill}
    \begin{subfigure}{0.48\textwidth}
        \centering
        \includegraphics[width=\linewidth]{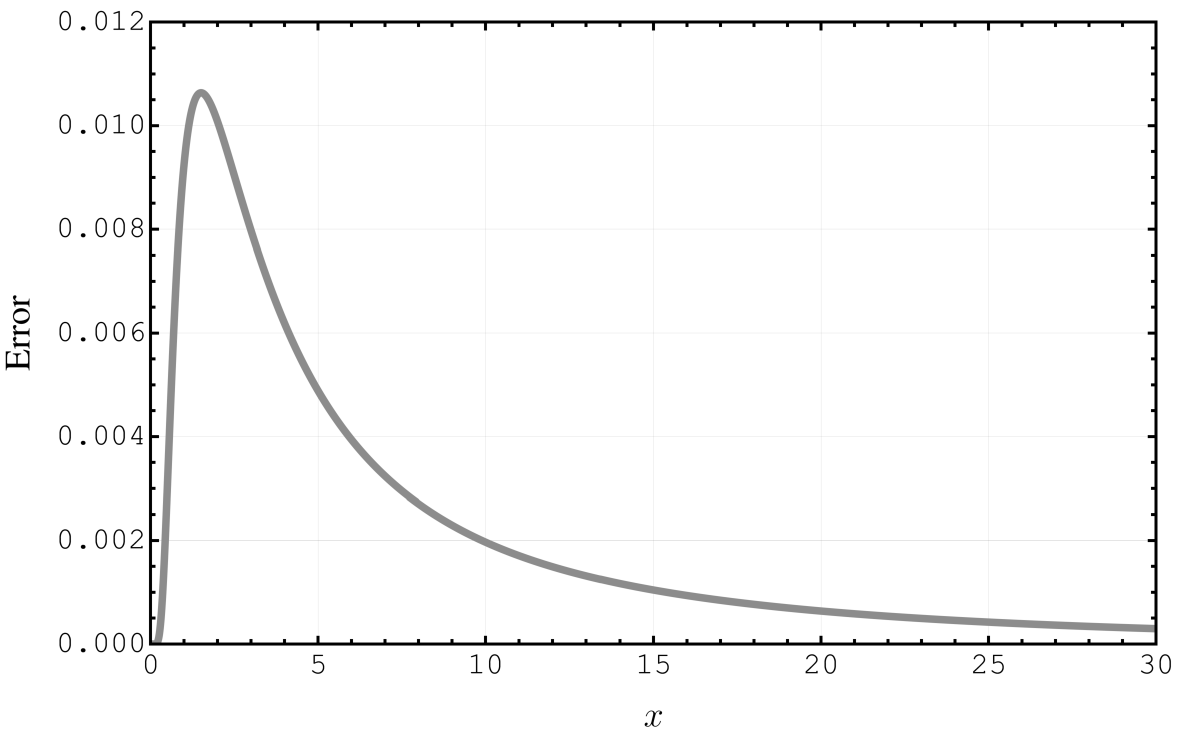}
        \caption{Corresponding lowerbound error.}
        \label{fig:QST-pdf-uprbnd1-err-A30}
    \end{subfigure}
    \caption{Quasi-stationary distribution's pdf, $q_{A}(x)$, its lowerbound $l_{A}^{(1)}(x)$, and the corresponding error---all as functions of $x\in[0,A]$ for $A=30$.}
    \label{fig:QST-pdf-lwrbnd1-perf-A30}
\end{figure}

Likewise, Figures~\ref{fig:QST-pdf-lwrbnd1-perf-A10},~\ref{fig:QST-pdf-lwrbnd1-perf-A20}, and~\ref{fig:QST-pdf-lwrbnd1-perf-A30} show $l_{A}^{(1)}(x)$ and $q_{A}(x)$ as functions of $x\in[0,A]$, also for $A=10$, $20$, and $30$, respectively. Again, we see that the discrepancy between $q_{A}(x)$ and $l_{A}^{(1)}(x)$ is fairly small, for all $x\in[0,A]$, and rapidly gets even smaller as $A$ increases.

We have seen numerical evidence that the lowerbound $l_{A}^{(1)}(x)$ and the upperbound $u_{A}^{(1)}(x)$ are both tight bounds. They are also tight in a more formal sense, namely in the sense that
\begin{align*}
\lim_{A\to+\infty}\left[l_{A}^{(1)}(x)\right]
&=
\lim_{A\to+\infty}\left[u_{A}^{(1)}(x)\right]
=
h(x),
\;
\text{pointwise, for each fixed $x\in\mathbb{R}$},
\end{align*}
where $h(x)$ is the stationary distribution's pdf~\eqref{eq:SR-StDist-def}--\eqref{eq:SR-StDist-answer}; here we used~\eqref{eq:Q-to-H-unif-conv} and~\eqref{eq:lambda-dbl-ineq}. Put another way, we see that $q_{A}(x)$ converges pointwise to $h(x)$, as $A\to+\infty$.




The pair $(l_{A}^{(1)}(x),u_{A}^{(1)}(x))$ can be used to lower- and upper-bound the cdf $Q_{A}(x)$. For example, inequality~\eqref{eq:QSD-pdf-qA-uprbnd1-def} may be rewritten as
\begin{align*}
\dfrac{q_{A}(x)}{Q_{A}(x)}
=
\dfrac{\partial}{\partial x}\log[Q_{A}(x)]
&\le
\dfrac{2}{x^2},
\;\;
\text{so that}
\;\;
\int_{x}^{A}\partial_t\log[Q_{A}(t)]
\le
2\int_{x}^{A}\dfrac{dt}{t^2},
\end{align*}
whence
\begin{align}\label{eq:QSD-cdf-lwrbnd1-def}
Q_{A}(x)
&\ge
e^{\tfrac{2}{A}}H(x)\eqqcolon L_{A}^{(1)}(x),
\;
x\in\mathbb{R},
\;
A>0,
\end{align}
which is a lowerbound for $Q_{A}(x)$ that was recently obtained by~\cite{Li+Polunchenko:SA2020} via a different argument. However, the argument used by~\cite{Li+Polunchenko:SA2020} to get the lowerbound $L_{A}^{(1)}(x)$ required that $A\ge\tilde{A}$. We can now see that that requirement can be lifted. Moreover, the bound $L_{A}^{(1)}(x)$ is a cdf in itself, with the corresponding density supported on $[0,A]$.



Likewise, the lowerbound $l_{A}^{(1)}(x)$ can be used to get a new upperbound for $Q_{A}(x)$. Specifically, from~\eqref{eq:QSD-pdf-qA-lwrbnd1-def} we obtain
\begin{align*}
\int_{x}^{A}\partial_t\log[Q_{A}(t)]
&\ge
 -\dfrac{2}{A}+\dfrac{2}{x}-2\log\left(\dfrac{A}{x}\right),
\end{align*}
whence
\begin{align}\label{eq:QSD-cdf-uprbnd1-def}
Q_{A}(x)
&\le
e^{\tfrac{2}{A}}H(x)\left(\dfrac{A}{x}\right)^{2\lambda}
\eqqcolon U_{A}^{(1)}(x),
\;
x\in[0,A],
\;
A>0.
\end{align}

Since for any $A\ge \tilde{A}$ we have $\xi\coloneqq\sqrt{1-8\lambda}\le1-4\lambda$, so that $1-\xi\ge 4\lambda$, we can conclude that the new upperbound $U_{A}^{(1)}(x)$ is tighter than the upperbound
\begin{align*}
Q_{A}(x)
&\le
e^{\tfrac{2}{A}}H(x)\left(\dfrac{A}{x}\right)^{\tfrac{1}{2}-\tfrac{1}{2}\xi},
\;\;
x\in[0,A],
\;
A\ge \tilde{A},
\end{align*}
previously obtained by~\cite{Li+Polunchenko:SA2020} via a different argument. Moreover, the new upperbound $U_{A}^{(1)}(x)$ is valid for {\em all} $A>0$, whereas the proof of the old upperbound due to~\cite{Li+Polunchenko:SA2020} does not go through unless $A\ge\tilde{A}$.


Now, let us try to get yet another representation for $q_{A}(x)$, different from~\eqref{eq:QSD-pdf-qA-via-cdf-QA-v1} and~\eqref{eq:QSD-pdf-qA-via-cdf-QA-v2}, but again with no other special functions involved other than the modified Bessel $K$ function. Such a representation will then allow us to appeal to some of the latest monotonicity properties of the modified Bessel $K$ function and functionals thereof, and establish new bounds for $q_{A}(x)$. First recall the identity
\begin{align*}
W_{\varkappa+1,b}(z)
&=
\left(\dfrac{z}{2}-\varkappa\right)W_{\varkappa,b}(z)-z\left[\dfrac{\partial}{\partial z}W_{\varkappa,b}(z)\right];
\end{align*}
cf., e.g.,~\cite[Identity~(2.4.24),~p.~25]{Slater:Book1960}. Hence
\begin{align*}
W_{1,b}(x)
&=
x\left\{\dfrac{1}{2}W_{0,b}(x)-\left[\dfrac{\partial}{\partial x}W_{0,b}(x)\right]\right\},
\end{align*}
which, on account of~\eqref{eq:Whit0-BesselK-id}, can subsequently be brought to the form
\begin{align}\label{eq:WhitW1b-BesselK-rel1}
W_{1,\tfrac{1}{2}\xi}\left(\dfrac{2}{x}\right)
&=
\dfrac{1}{x}\sqrt{\dfrac{2}{\pi x}}\left\{\left(1-\dfrac{x}{2}\right)K_{\tfrac{1}{2}\xi}\left(\dfrac{1}{x}\right)-\left.\left[\dfrac{\partial}{\partial u}K_{\tfrac{1}{2}\xi}(u)\right]\right|_{u=\tfrac{1}{x}}\right\}.
\end{align}

To proceed, we appeal to the derivative formulae
\begin{align}\label{eq:BesselK-deriv}
\dfrac{\partial}{\partial z}K_{b}(z)
&=
-K_{b-1}(z)-\dfrac{b}{z}K_{b}(z)
=
-K_{b+1}(z)+\dfrac{b}{z}K_{b}(z),
\end{align}
given, e.g., by~\cite[Identities~8.486.12~and~8.486.13,~p.~938]{Gradshteyn+Ryzhik:Book2014}. By substituting~\eqref{eq:BesselK-deriv} back over into~\eqref{eq:WhitW1b-BesselK-rel1} we arrive at two additional equivalent expressions for $W_{1,\xi/2}(2/x)$ in terms of the Bessel $K$ function:
\begin{align}
W_{1,\tfrac{1}{2}\xi}\left(\dfrac{2}{x}\right)
&=
\dfrac{1}{x}\sqrt{\dfrac{2}{\pi x}}\left\{\left[1-\dfrac{x}{2}\left(1-\xi\right)\right]K_{\tfrac{1}{2}\xi}\left(\dfrac{1}{x}\right)+K_{\tfrac{1}{2}\xi-1}\left(\dfrac{1}{x}\right)\right\},
\label{eq:WhitW1b-BesselK-rel2}
\\
W_{1,\tfrac{1}{2}\xi}\left(\dfrac{2}{x}\right)
&=
\dfrac{1}{x}\sqrt{\dfrac{2}{\pi x}}\left\{\left[1-\dfrac{x}{2}\left(1+\xi\right)\right]K_{\tfrac{1}{2}\xi}\left(\dfrac{1}{x}\right)+K_{\tfrac{1}{2}\xi+1}\left(\dfrac{1}{x}\right)\right\}.
\label{eq:WhitW1b-BesselK-rel3}
\end{align}

We hasten to note the identities
\begin{align}
\left.\left[\dfrac{\partial}{\partial u}K_{\tfrac{1}{2}\xi}(u)\right]\right|_{u=\tfrac{1}{A}}
&=
\left[1-\dfrac{A}{2}\right]K_{\tfrac{1}{2}\xi}\left(\dfrac{1}{A}\right),
\nonumber
\\
K_{\tfrac{1}{2}\xi-1}\left(\dfrac{1}{A}\right)
&=
\left[\dfrac{A}{2}(1-\xi)-1\right]K_{\tfrac{1}{2}\xi}\left(\dfrac{1}{A}\right),
\label{eq:qA-zero-BesselK-ind1down-id}
\\
K_{\tfrac{1}{2}\xi+1}\left(\dfrac{1}{A}\right)
&=
\left[\dfrac{A}{2}(1+\xi)-1\right]K_{\tfrac{1}{2}\xi}\left(\dfrac{1}{A}\right),
\label{eq:qA-zero-BesselK-ind1up-id}
\end{align}
which follow at once from~\eqref{eq:WhitW1b-BesselK-rel1},~\eqref{eq:WhitW1b-BesselK-rel2}, and~\eqref{eq:WhitW1b-BesselK-rel3}, respectively, on account of~\eqref{eq:lambda-eqn}; all three hold for {\em any} $A>0$. Also, identity~\eqref{eq:qA-zero-BesselK-ind1down-id} and identity~\eqref{eq:qA-zero-BesselK-ind1up-id} are the effectively equivalent to one another, because the Bessel $K$ function is an even function of its index.

If we now substitute~\eqref{eq:WhitW1b-BesselK-rel1},~\eqref{eq:WhitW1b-BesselK-rel2}, and~\eqref{eq:WhitW1b-BesselK-rel3} back into~\eqref{eq:QSD-pdf-answer} we will get the following three equivalent expressions for $q_{A}(x)$, all in terms of the modified Bessel $K$ function:
\begin{align*}
q_{A}(x)
&=
\dfrac{1}{x^2}\left\{\left[1-\dfrac{x}{2}\right]Q_{A}(x)-\dfrac{e^{-\tfrac{1}{x}}\,\sqrt{\dfrac{2}{\pi x}}\,\left.\left[\dfrac{\partial}{\partial u}K_{\tfrac{1}{2}\xi}(u)\right]\right|_{u=\tfrac{1}{x}}}{e^{-\tfrac{1}{A}}\,\sqrt{\dfrac{2}{\pi A}}\,K_{\tfrac{1}{2}\xi}\left(\dfrac{1}{A}\right)}\right\}\indicator{x\in[0,A]}
\\
&=
\dfrac{1}{x^2}\left\{\left[1-\dfrac{x}{2}(1-\xi)\right]Q_{A}(x)+\dfrac{e^{-\tfrac{1}{x}}\,\sqrt{\dfrac{2}{\pi x}}\,K_{\tfrac{1}{2}\xi-1}\left(\dfrac{1}{x}\right)}{e^{-\tfrac{1}{A}}\,\sqrt{\dfrac{2}{\pi A}}\,K_{\tfrac{1}{2}\xi}\left(\dfrac{1}{A}\right)}\right\}\indicator{x\in[0,A]}
\\
&=
\dfrac{1}{x^2}\left\{\left[1-\dfrac{x}{2}(1+\xi)\right]Q_{A}(x)+\dfrac{e^{-\tfrac{1}{x}}\,\sqrt{\dfrac{2}{\pi x}}\,K_{\tfrac{1}{2}\xi+1}\left(\dfrac{1}{x}\right)}{e^{-\tfrac{1}{A}}\,\sqrt{\dfrac{2}{\pi A}}\,K_{\tfrac{1}{2}\xi}\left(\dfrac{1}{A}\right)}\right\}\indicator{x\in[0,A]}.
\end{align*}

The plan now is to use the foregoing three expressions for $q_{A}(x)$ in conjunction with certain recently established monotonicity properties of the modified Bessel $K$ function and its functions to get new bounds for $q_{A}(x)$. To that end, consider first~\cite[Proposition~4.5,~p.~2957]{Yang+Zheng:PAMS2017} whereby the function
\begin{align*}
f_1(x;b)
&\coloneqq
\dfrac{x}{1+x}\dfrac{K_{b+1}(x)}{K_{b}(x)}
\end{align*}
is strictly increasing in $x$ on $(0,+\infty)$ for any $b\in(0,1/2)$. This gives
\begin{align*}
\dfrac{1}{1+x}\left[K_{\tfrac{1}{2}\xi+1}\left(\dfrac{1}{x}\right)\left/K_{\tfrac{1}{2}\xi}\left(\dfrac{1}{x}\right)\right.\right]
&>
\dfrac{1}{1+A}\left[K_{\tfrac{1}{2}\xi+1}\left(\dfrac{1}{A}\right)\left/K_{\tfrac{1}{2}\xi}\left(\dfrac{1}{A}\right)\right.\right], \\[1ex]
& \qquad x\in(0,A),\;\; A\ge\tilde{A},
\end{align*}
because $\xi\in[0,1]$ for $A\ge\tilde{A}$.

After some elementary algebra this translates to
\begin{align*} 
l_{A}^{(3)}(x)
&\coloneqq
\dfrac{2}{x^2} Q_{A}(x)\dfrac{3+\xi}{4}\dfrac{A-x}{A+1}\indicator{x\in[0,A]}
\le
q_{A}(x),
\;
x\in\mathbb{R},
\;
A\ge\tilde{A}.
\end{align*}

\begin{figure}[h!]
    \centering
    \begin{subfigure}{0.48\textwidth}
        \centering
        \includegraphics[width=\linewidth]{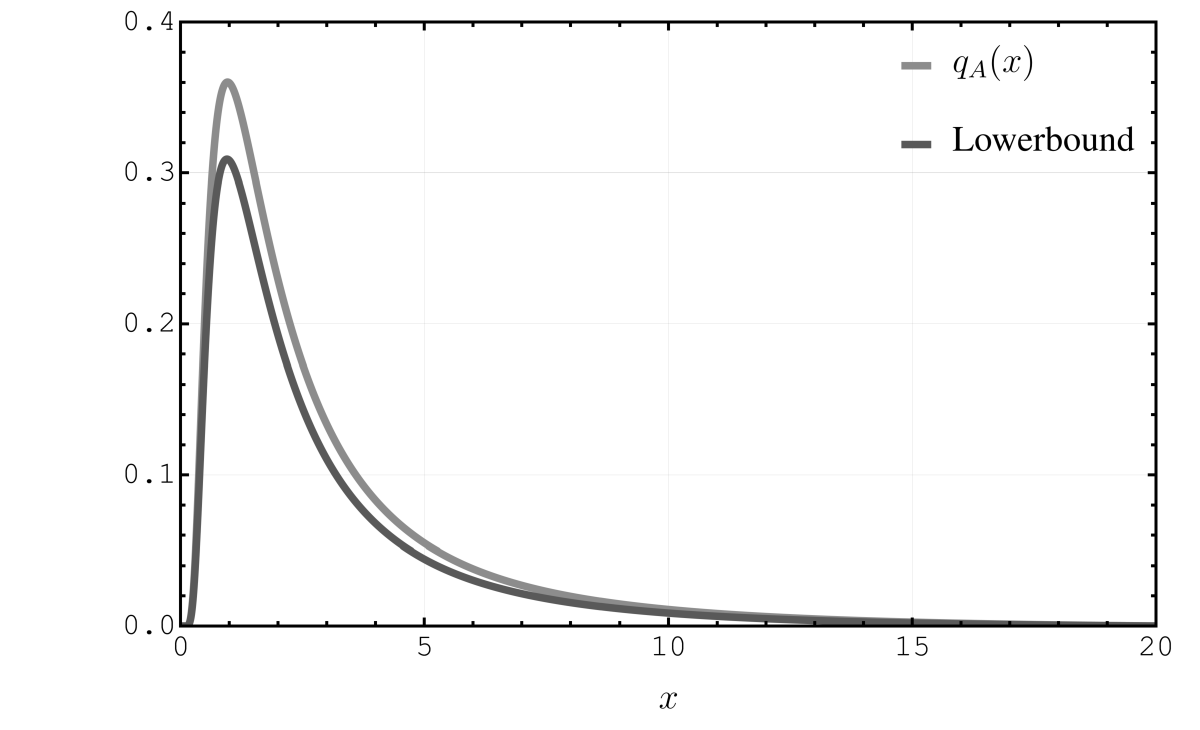}
        \caption{$q_{A}(x)$ and $l_{A}^{(3)}(x)$.}
        \label{fig:QST-pdf-lwrbnd3-A20}
    \end{subfigure}
    \hspace*{\fill}
    \begin{subfigure}{0.48\textwidth}
        \centering
        \includegraphics[width=\linewidth]{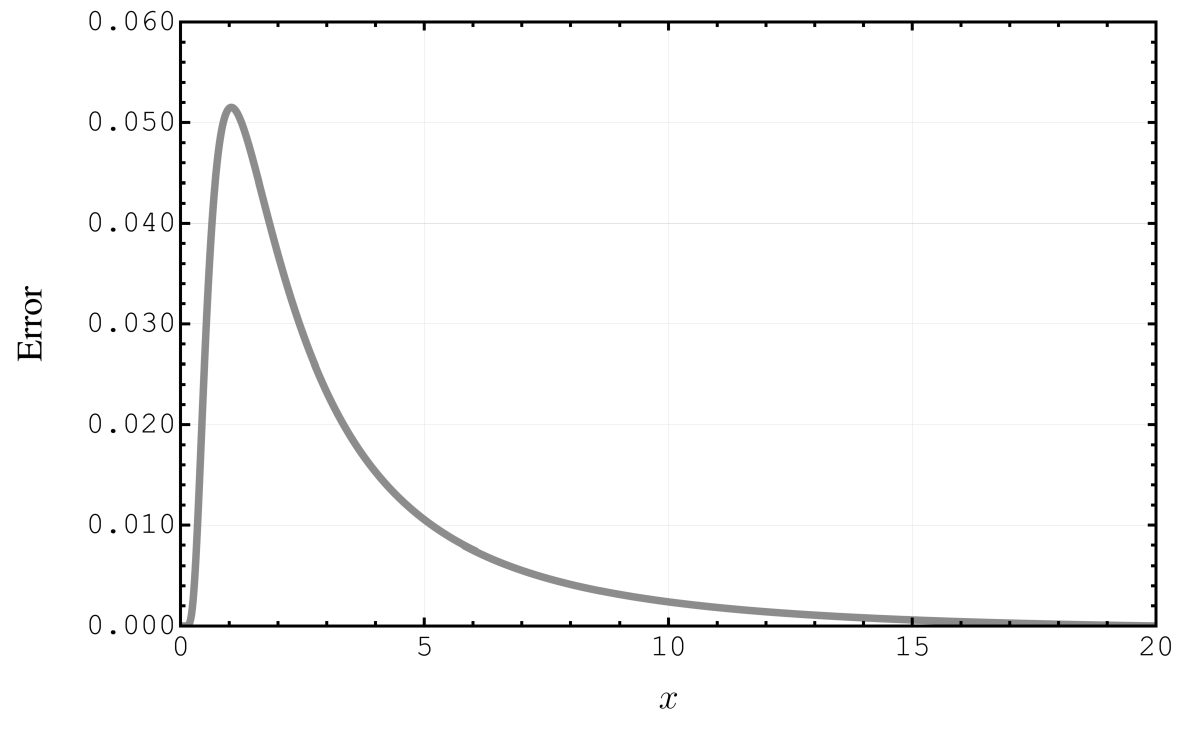}
        \caption{Corresponding lowerbound error.}
        \label{fig:QST-pdf-uprbnd3-err-A20}
    \end{subfigure}
    \caption{Quasi-stationary distribution's pdf, $q_{A}(x)$, its lowerbound $l_{A}^{(3)}(x)$, and the corresponding error---all as functions of $x\in[0,A]$ for $A=20$.}
    \label{fig:QST-pdf-lwrbnd3-perf-A20}
\end{figure}

\begin{figure}[h!]
    \centering
    \begin{subfigure}{0.48\textwidth}
        \centering
        \includegraphics[width=\linewidth]{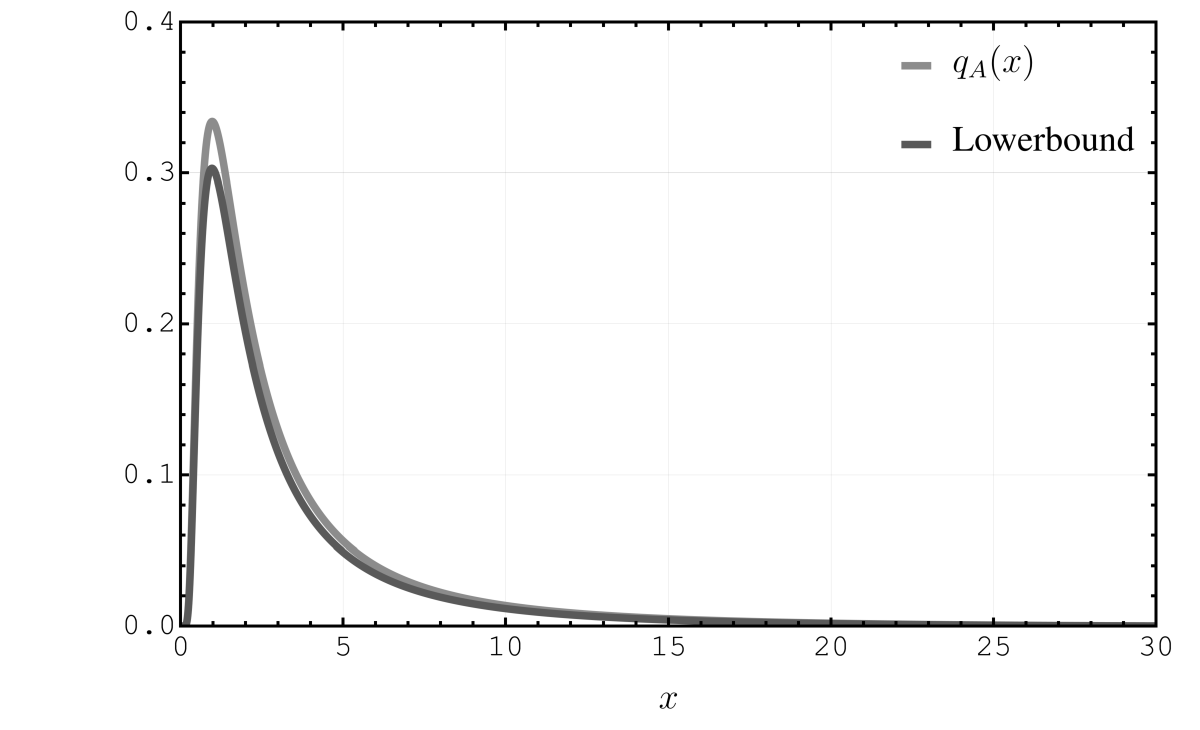}
        \caption{$q_{A}(x)$ and $l_{A}^{(3)}(x)$.}
        \label{fig:QST-pdf-lwrbnd3-A30}
    \end{subfigure}
    \hspace*{\fill}
    \begin{subfigure}{0.48\textwidth}
        \centering
        \includegraphics[width=\linewidth]{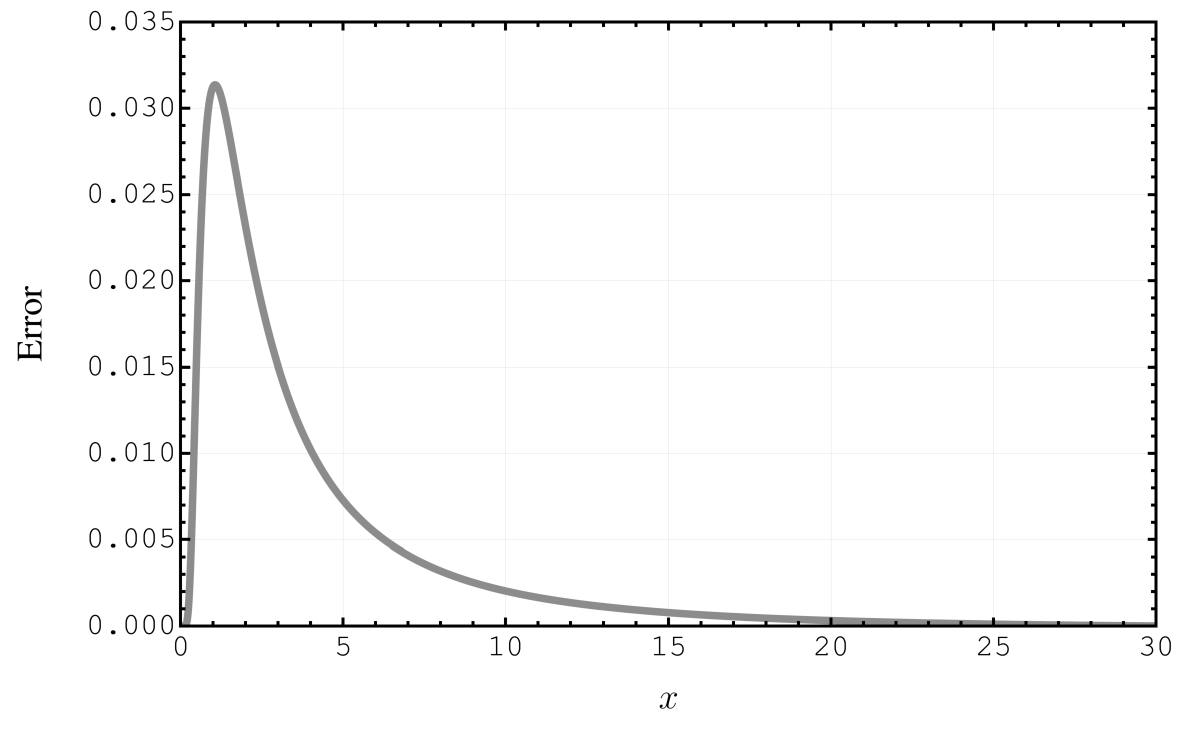}
        \caption{Corresponding lowerbound error.}
        \label{fig:QST-pdf-uprbnd3-err-A30}
    \end{subfigure}
    \caption{Quasi-stationary distribution's pdf, $q_{A}(x)$, its lowerbound $l_{A}^{(3)}(x)$, and the corresponding error---all as functions of $x\in[0,A]$ for $A=30$.}
    \label{fig:QST-pdf-lwrbnd3-perf-A30}
\end{figure}

Figures~\ref{fig:QST-pdf-lwrbnd3-A20} and~\ref{fig:QST-pdf-lwrbnd3-A30} show the performance of $l_{A}^{(3)}(x)$ as a function of $x\in[0,A]$ for $A=20$ and $A=30$, respectively; note that $l_{A}^{(3)}(x)$ requires $A\ge\tilde{A}$. We see that the bound is generally looser than $l_{A}^{(1)}(x)$, unless $x$ is sufficiently close to $A$. We also see that the tightness of the bound improves as $A$ increases.

The lowerbound $l_{A}^{(3)}(x)$ vanishes at $x=0$ as well as at $x=A$, just as the pdf $q_{A}(x)$ itself does. Moreover, the bound is also tight in the sense that
\begin{align*}
\lim_{A\to+\infty}\left[l_{A}^{(3)}(x)\right]
&=
h(x),
\;
\text{pointwise, for any {\em fixed}}
\;
x\ge0,
\end{align*}
because $\lim_{A\to+\infty}\xi=1$ and $\sup_{x\in\mathbb{R}}|Q_{A}(x)-H(x)|\to 0$, as $A\to+\infty$.

Now, consider~\cite[Lemma~2.2,~p.~376]{Gaunt:JMAA2014} which states, in particular, that the function
\begin{align*}
f_2(x;b)
&\coloneqq
\dfrac{K_{b-1}(x)}{K_{b}(x)}
\end{align*}
is strictly monotone decreasing in $x$ on $(0,+\infty)$ for every fixed $b<1/2$. This gives
\begin{align*}
\left[K_{\tfrac{1}{2}\xi-1}\left(\dfrac{1}{x}\right)\left/K_{\tfrac{1}{2}\xi}\left(\dfrac{1}{x}\right)\right.\right]
&>
\left[K_{\tfrac{1}{2}\xi-1}\left(\dfrac{1}{A}\right)\left/K_{\tfrac{1}{2}\xi}\left(\dfrac{1}{A}\right)\right.\right],
\;\;
x\in(0,A],
\;
A\ge\tilde{A},
\end{align*}
whence, in view of~\eqref{eq:qA-zero-BesselK-ind1up-id}, we find
\begin{align*}
\left[\dfrac{A}{2}(1-\xi)-1\right]K_{\tfrac{1}{2}\xi}\left(\dfrac{1}{x}\right)
&<
K_{\tfrac{1}{2}\xi-1}\left(\dfrac{1}{x}\right),
\;\;
x\in(0,A],
\;
A\ge\tilde{A},
\end{align*}
and consequently
\begin{align*} 
q_{A}(x)
&\le
\dfrac{2}{x^2}Q_{A}(x)\dfrac{1-\xi}{2}\left(A-x\right)\indicator{x\in[0,A]}
\eqqcolon u_{A}^{(3)}(x),
\;
x\in\mathbb{R},
\;
A\ge\tilde{A}.
\end{align*}

\begin{figure}[h!]
    \centering
    \begin{subfigure}{0.48\textwidth}
        \centering
        \includegraphics[width=\linewidth]{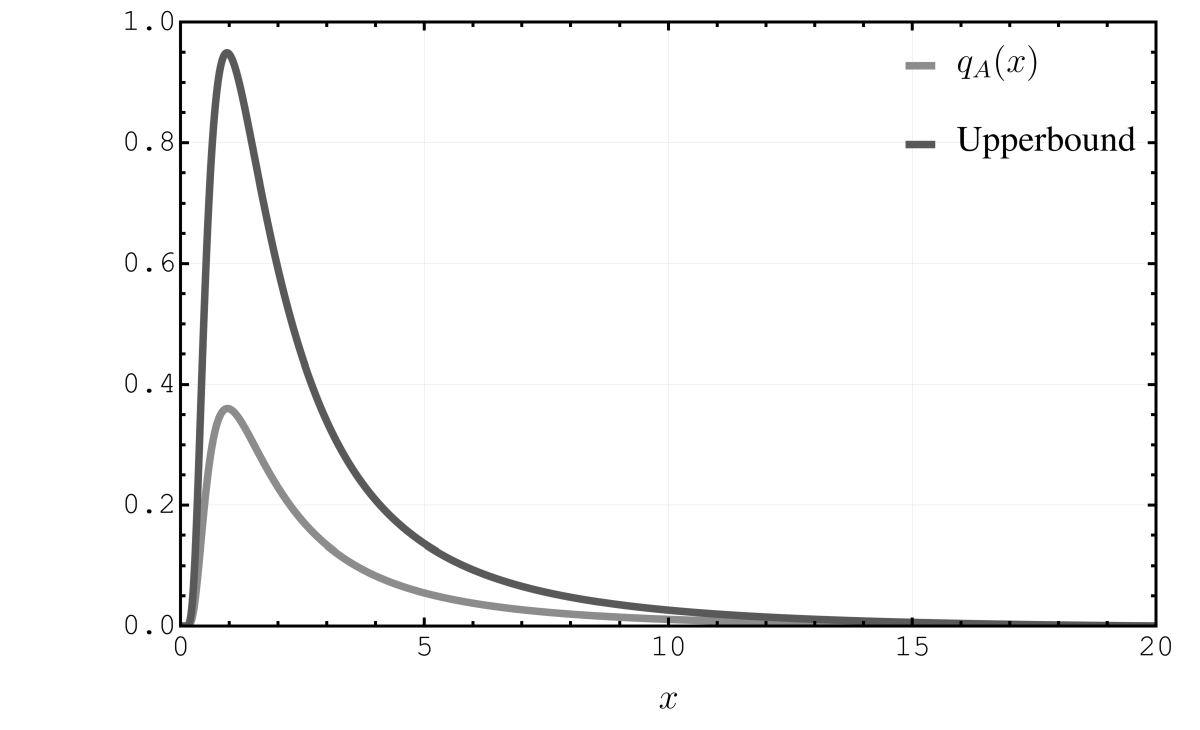}
        \caption{$q_{A}(x)$ and $u_{A}^{(3)}(x)$.}
        \label{fig:QST-pdf-uprbnd3-A20}
    \end{subfigure}
    \hspace*{\fill}
    \begin{subfigure}{0.48\textwidth}
        \centering
        \includegraphics[width=\linewidth]{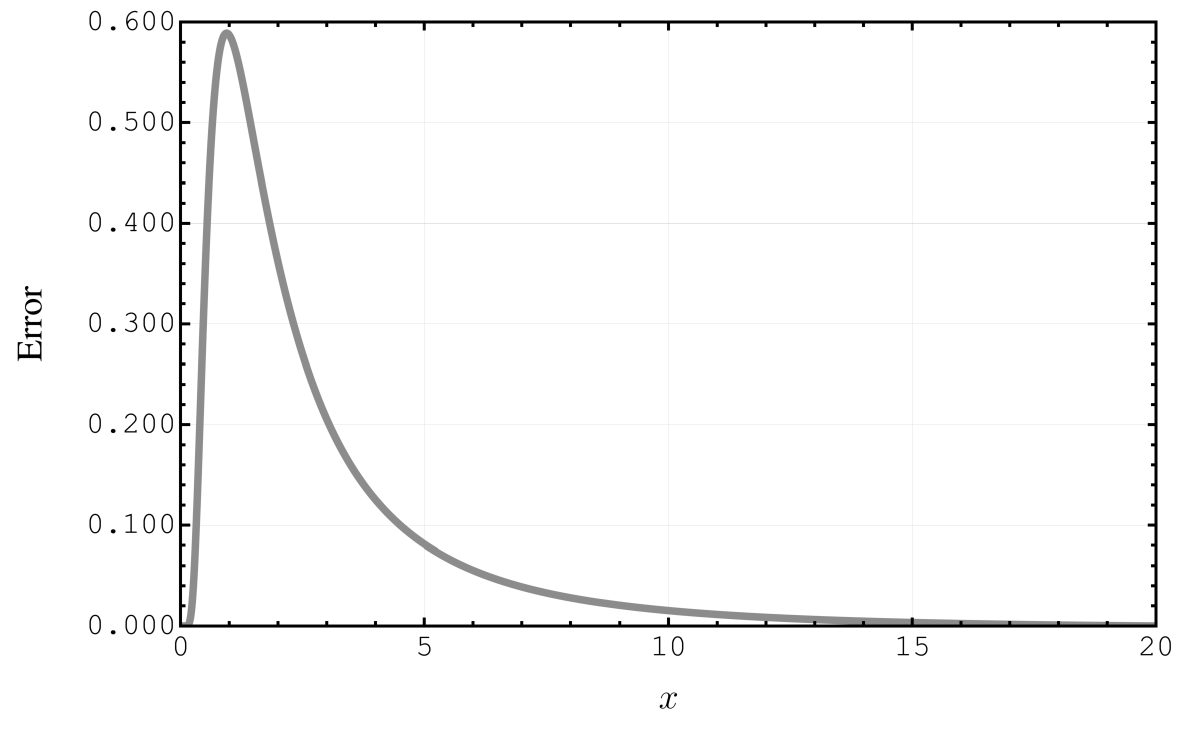}
        \caption{Corresponding upperbound error.}
        \label{fig:QST-pdf-uprbnd3-err-A20}
    \end{subfigure}
    \caption{Quasi-stationary distribution's pdf, $q_{A}(x)$, its upperbound $u_{A}^{(3)}(x)$, and the corresponding error---all as functions of $x\in[0,A]$ for $A=20$.}
    \label{fig:QST-pdf-uprbnd3-perf-A20}
\end{figure}

\begin{figure}[h!]
    \centering
    \begin{subfigure}{0.48\textwidth}
        \centering
        \includegraphics[width=\linewidth]{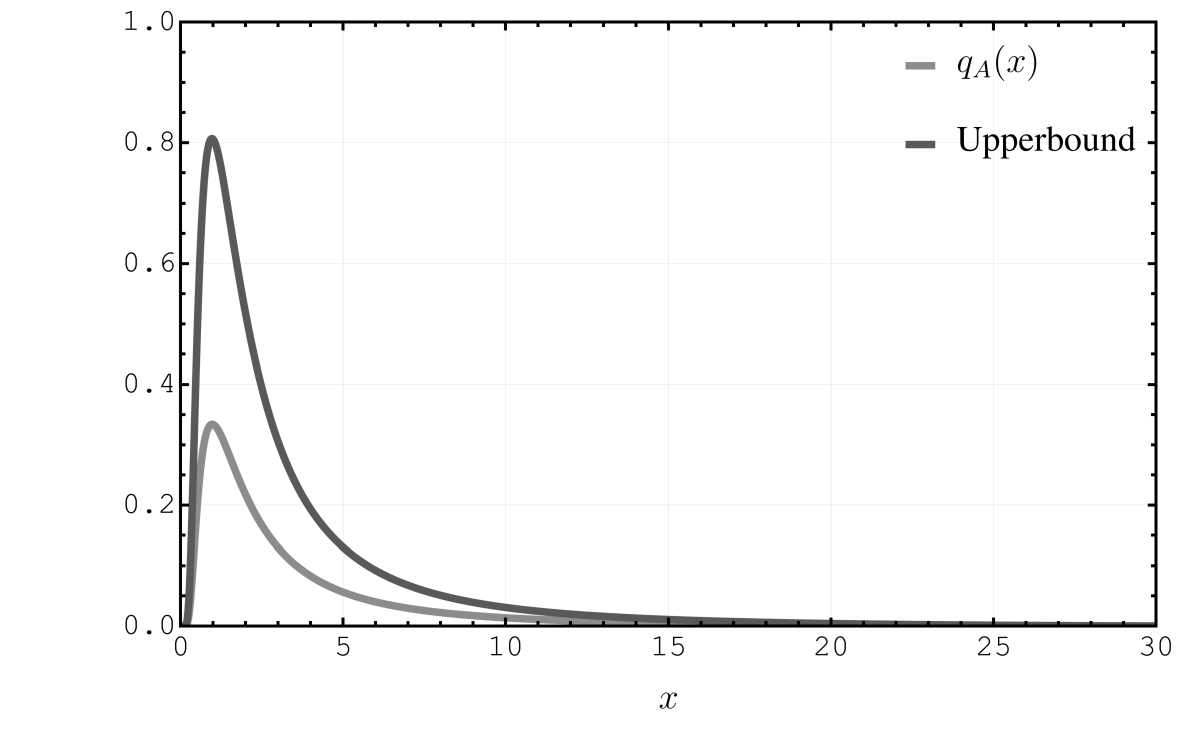}
        \caption{$q_{A}(x)$ and $u_{A}^{(3)}(x)$.}
        \label{fig:QST-pdf-uprbnd3-A30}
    \end{subfigure}
    \hspace*{\fill}
    \begin{subfigure}{0.48\textwidth}
        \centering
        \includegraphics[width=\linewidth]{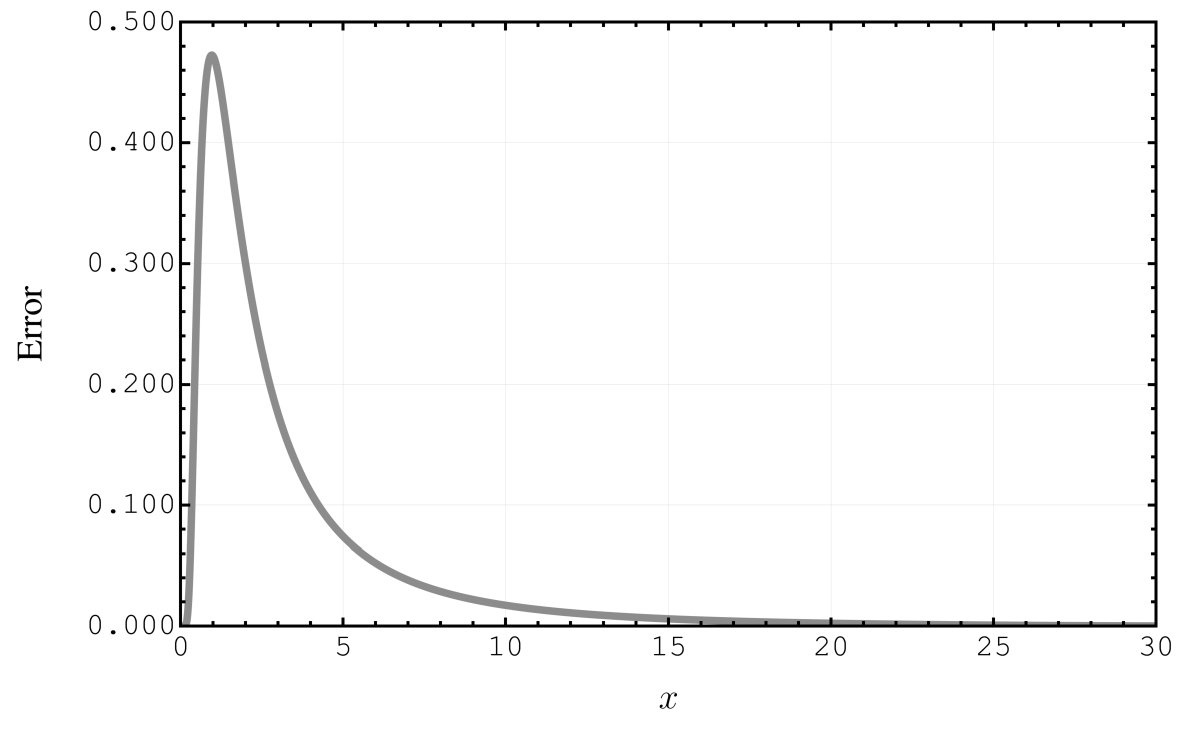}
        \caption{Corresponding upperbound error.}
        \label{fig:QST-pdf-uprbnd3-err-A30}
    \end{subfigure}
    \caption{Quasi-stationary distribution's pdf, $q_{A}(x)$, its upperbound $u_{A}^{(3)}(x)$, and the corresponding error---all as functions of $x\in[0,A]$ for $A=30$.}
    \label{fig:QST-pdf-uprbnd3-perf-A30}
\end{figure}

The upperbound $u_{A}^{(3)}(x)$ works only for $A\ge \tilde{A}$, while the upperbound $u_{A}^{(1)}(x)$ given by~\eqref{eq:QSD-pdf-qA-uprbnd1-def} works for any $A>0$. The upperbound $u_{A}^{(3)}(x)$ has the property that $u_{A}^{(3)}(0)=u_{A}^{(3)}(A)=0$, which the upperbound $u_{A}^{(1)}(x)$ given by~\eqref{eq:QSD-pdf-qA-uprbnd1-def} lacks. However, the bound $u_{A}^{(3)}(x)$ is tighter than $u_{A}^{(1)}(x)$ only for values of $x$ that sufficiently close to $A$; otherwise $u_{A}^{(3)}(x)$ is looser than $u_{A}^{(1)}(x)$ for values of $x$ that are sufficiently far to the left of $A$. Figures~\ref{fig:QST-pdf-uprbnd3-perf-A20} and~\ref{fig:QST-pdf-uprbnd3-perf-A30} show the performance of $u_{A}^{(3)}(x)$ for $A=20$ and $30$. We also have $\lim_{A\to+\infty}\left[u_{A}^{(3)}(x)\right]=h(x)$ for any fixed $x\ge0$. This is because $\lim_{A\to+\infty}(1-\xi)=0$ but $\lim_{A\to+\infty}[(1-\xi)A]=2$, as can be seen from~\eqref{eq:lambda-dbl-ineq}.

%


Next, from~\cite[Inequality~(4.7),~p.~2956]{Yang+Zheng:PAMS2017}, whereby
\begin{align*}
\dfrac{1}{x}\dfrac{K_{b-1}(x)}{K_{b}(x)}
&>
\dfrac{1}{x+1}\dfrac{K_{b+1}(x)}{K_{b}(x)},
\;\;
x>0,
\;\;
b\in(0,1/2),
\end{align*}
we find that
\begin{align*}
\dfrac{e^{-\tfrac{1}{x}}\,\sqrt{\dfrac{2}{\pi x}}\,K_{\tfrac{1}{2}\xi-1}\left(\dfrac{1}{x}\right)}{e^{-\tfrac{1}{A}}\,\sqrt{\dfrac{2}{\pi A}}\,K_{\tfrac{1}{2}\xi}\left(\dfrac{1}{A}\right)}
>
\xi Q_{A}(x),
\end{align*}
and therefore
\begin{align*} 
l_{A}^{(4)}(x)
&\coloneqq
\dfrac{Q_{A}(x)}{x^2}\left[1+\xi-\dfrac{x}{2}(1-\xi)\right]\indicator{x\in[0,A]}
\le
q_{A}(x),
\;
x\in\mathbb{R},
\;
A>\tilde{A}.
\end{align*}

\begin{figure}[h!]
    \centering
    \begin{subfigure}{0.48\textwidth}
        \centering
        \includegraphics[width=\linewidth]{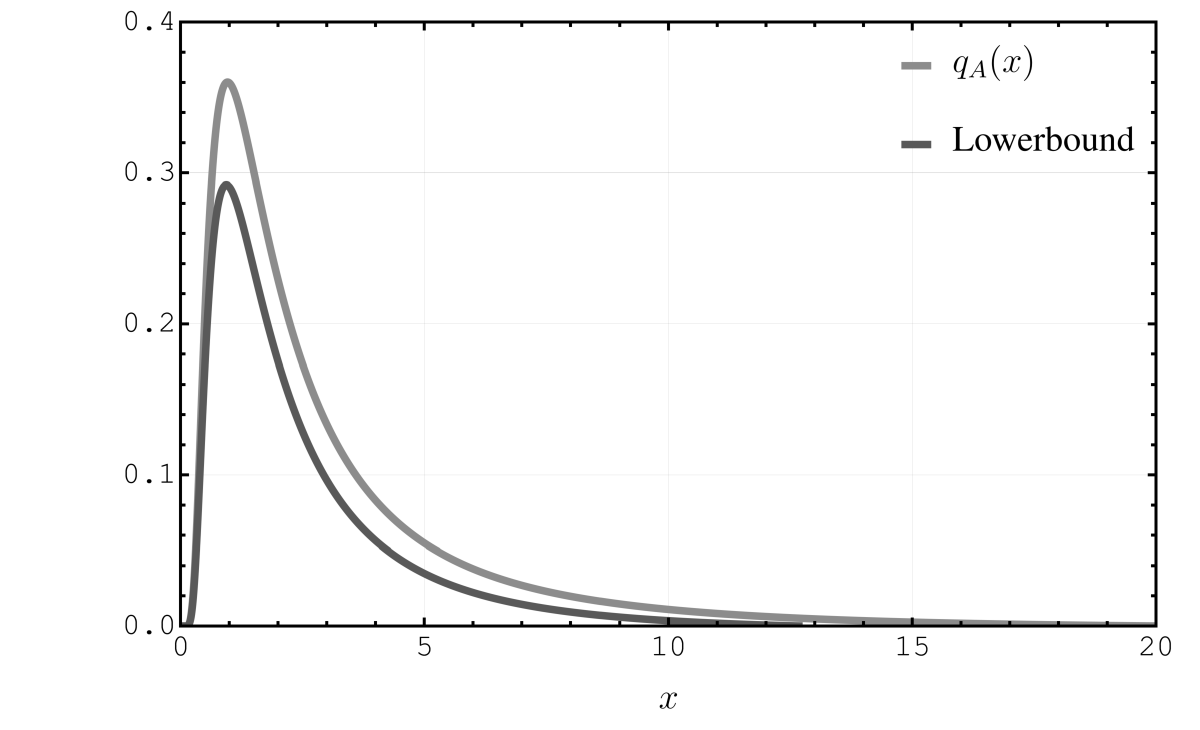}
        \caption{$q_{A}(x)$ and $l_{A}^{(4)}(x)$.}
        \label{fig:QST-pdf-lwrbnd4-A20}
    \end{subfigure}
    \hspace*{\fill}
    \begin{subfigure}{0.48\textwidth}
        \centering
        \includegraphics[width=\linewidth]{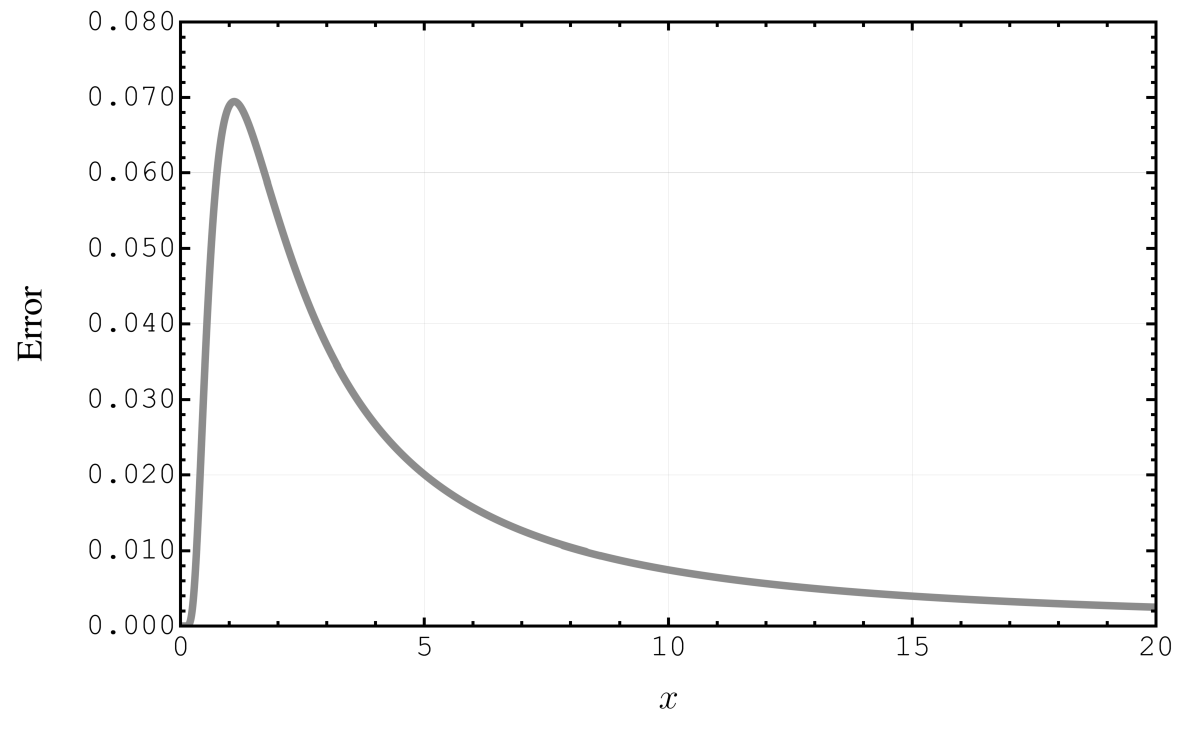}
        \caption{Corresponding lowerbound error.}
        \label{fig:QST-pdf-uprbnd4-err-A20}
    \end{subfigure}
    \caption{Quasi-stationary distribution's pdf, $q_{A}(x)$, its lowerbound $l_{A}^{(4)}(x)$, and the corresponding error---all as functions of $x\in[0,A]$ for $A=20$.}
    \label{fig:QST-pdf-lwrbnd4-perf-A20}
\end{figure}

\begin{figure}[h!]
    \centering
    \begin{subfigure}{0.48\textwidth}
        \centering
        \includegraphics[width=\linewidth]{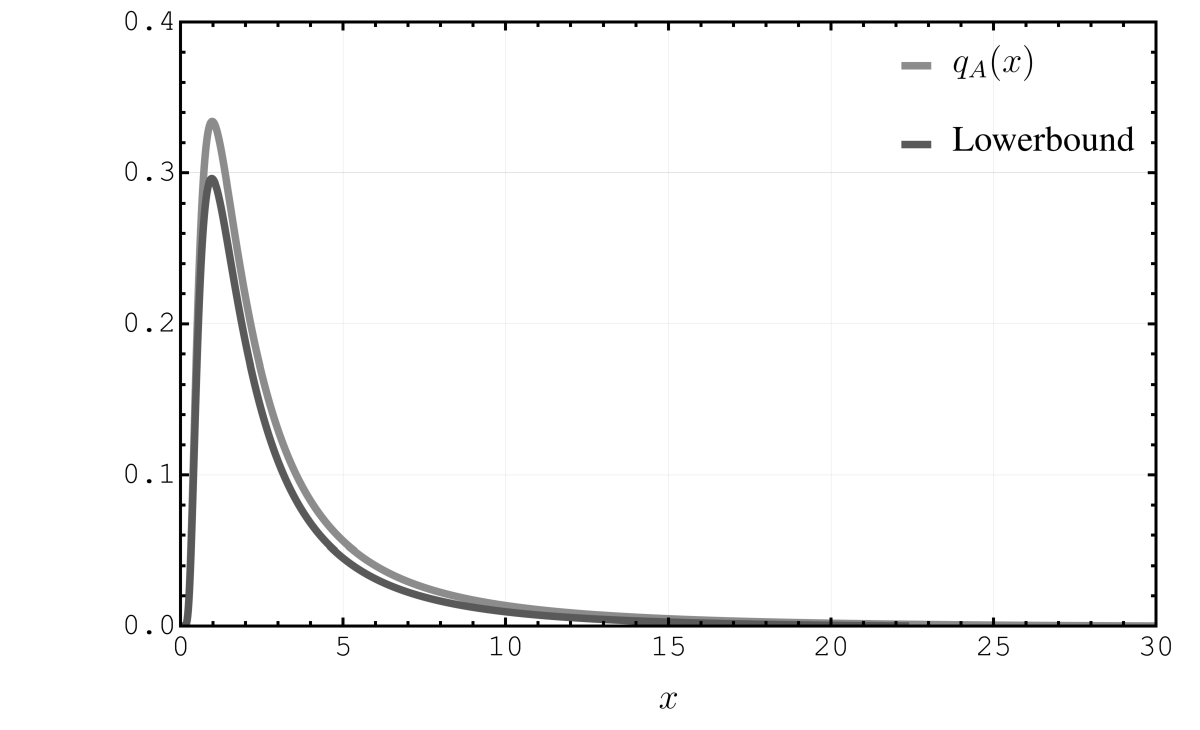}
        \caption{$q_{A}(x)$ and $l_{A}^{(4)}(x)$.}
        \label{fig:QST-pdf-lwrbnd4-A30}
    \end{subfigure}
    \hspace*{\fill}
    \begin{subfigure}{0.48\textwidth}
        \centering
        \includegraphics[width=\linewidth]{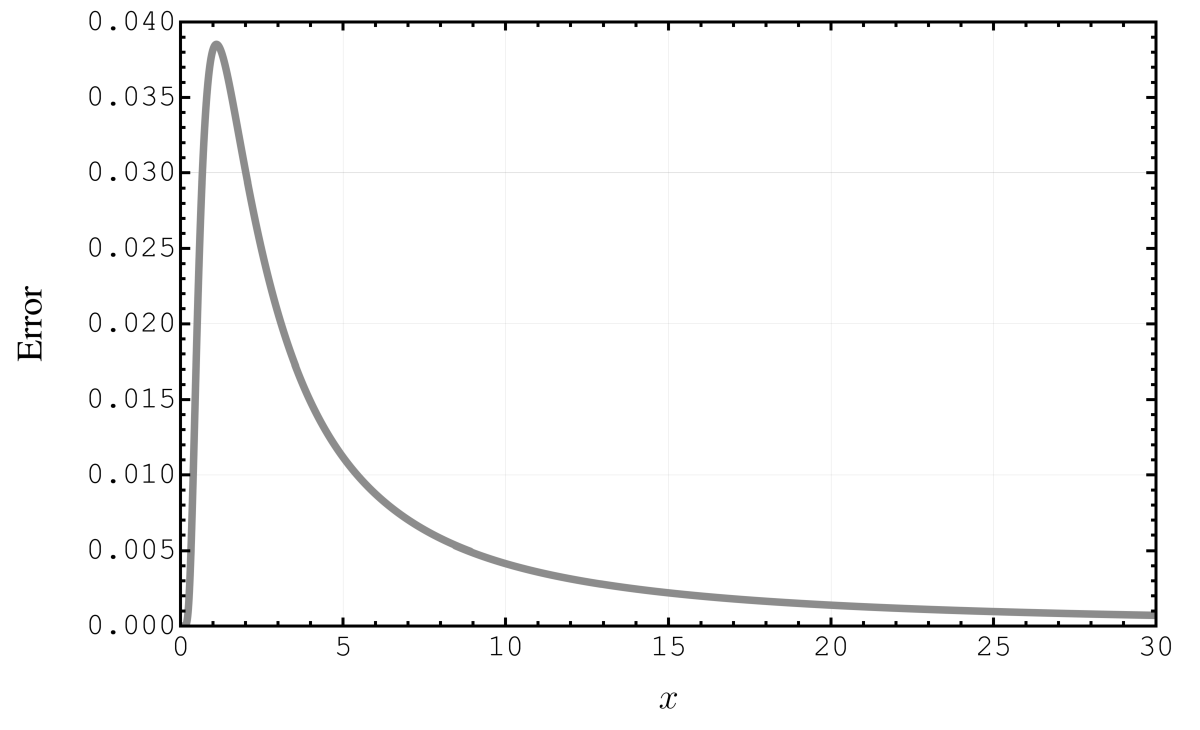}
        \caption{Corresponding lowerbound error.}
        \label{fig:QST-pdf-uprbnd4-err-A30}
    \end{subfigure}
    \caption{Quasi-stationary distribution's pdf, $q_{A}(x)$, its lowerbound $l_{A}^{(4)}(x)$, and the corresponding error---all as functions of $x\in[0,A]$ for $A=30$.}
    \label{fig:QST-pdf-lwrbnd4-perf-A30}
\end{figure}

Figures~\ref{fig:QST-pdf-lwrbnd4-perf-A20} and~\ref{fig:QST-pdf-lwrbnd4-perf-A30} show the performance of $l_{A}^{(4)}(x)$ as a function of $x\in[0,A]$ for $A=20$ and $A=30$. The figures suggest that $l_{A}^{(4)}(x)$ is an adequate lowerbound.

To proceed, we turn to~\cite[Property~4,~p.~2999]{Yang+Zheng:PAMS2022} which states that the function
\begin{align*}
f_3(x;b)
&\coloneqq
x\dfrac{K_{b-1}(x)}{K_{b}(x)}-x
\end{align*}
is increasing in $x$ on $(0,+\infty)$ for $b\in(-1/2,1/2)$. This gives
\begin{align*}
  \dfrac{1}{x}\left[K_{\tfrac{1}{2}\xi-1}\left(\dfrac{1}{x}\right)\left/K_{\tfrac{1}{2}\xi}\left(\dfrac{1}{x}\right)\right.\right]
  - \dfrac{1}{x}
  & \ge
  \dfrac{1}{A}\left[K_{\tfrac{1}{2}\xi-1}\left(\dfrac{1}{A}\right)\left/K_{\tfrac{1}{2}\xi}\left(\dfrac{1}{A}\right)\right.\right]
  -
  \dfrac{1}{A}, \\[1ex]
  & x\in(0,A], \; A\ge\tilde{A},
\end{align*}
which, in view of~\eqref{eq:qA-zero-BesselK-ind1down-id}, becomes
\begin{align*}
K_{\tfrac{1}{2}\xi-1}\left(\dfrac{1}{x}\right)
&\ge
\left[1-2\dfrac{x}{A}+\dfrac{x}{2}(1-\xi)\right]K_{\tfrac{1}{2}\xi}\left(\dfrac{1}{x}\right),
\;\;
x\in(0,A],
\;
A\ge\tilde{A}.
\end{align*}

Thus we obtain:
\begin{align}\label{eq:QSD-pdf-qA-lwrbnd5-def}
l_{A}^{(5)}(x)
&\coloneqq
\dfrac{2}{x^2} Q_{A}(x)\left(1-\dfrac{x}{A}\right)\indicator{x\in[0,A]}
\le
q_{A}(x),
\;
x\in\mathbb{R},
\;
A\ge\tilde{A},
\end{align}
which has the desired properties that $l_{A}^{(5)}(0)=l_{A}^{(5)}(A)=0$, and is a tight bound in the sense that $\lim_{A\to+\infty}l_{A}^{(5)}(x)=h(x)$ for any fixed $x\ge0$.

\begin{figure}[h!]
    \centering
    \begin{subfigure}{0.48\textwidth}
        \centering
        \includegraphics[width=\linewidth]{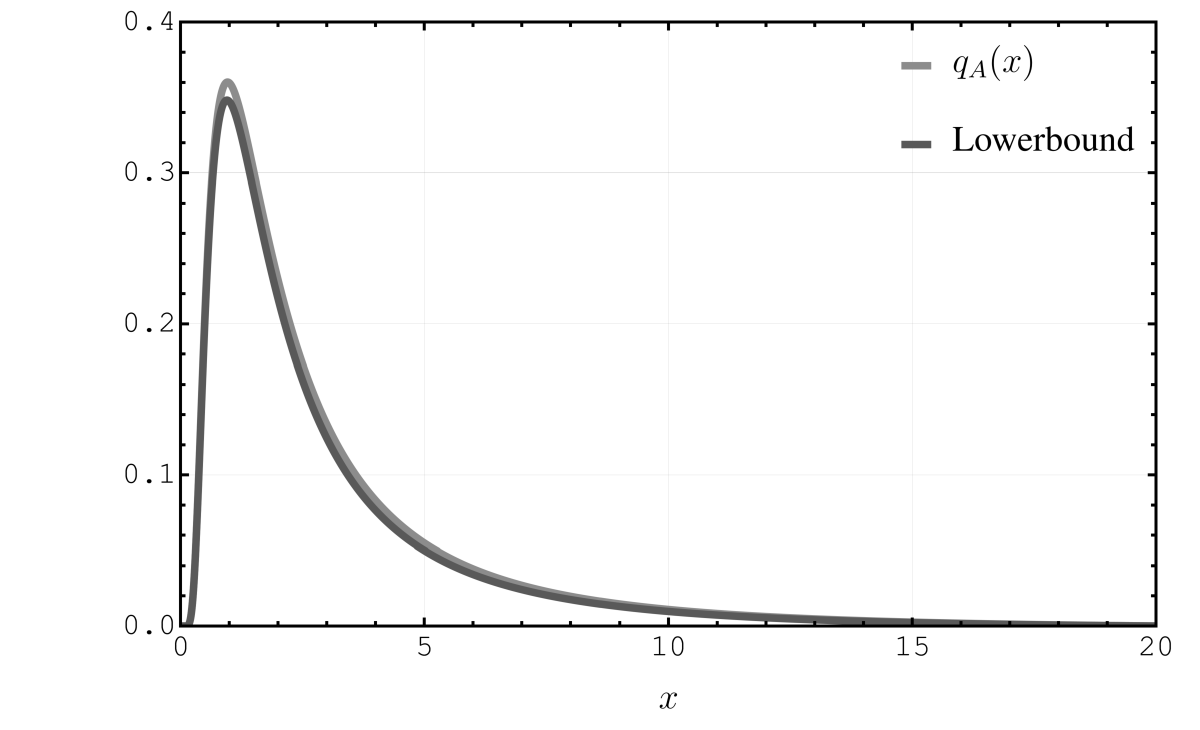}
        \caption{$q_{A}(x)$ and $l_{A}^{(5)}(x)$.}
        \label{fig:QST-pdf-lwrbnd5-A20}
    \end{subfigure}
    \hspace*{\fill}
    \begin{subfigure}{0.48\textwidth}
        \centering
        \includegraphics[width=\linewidth]{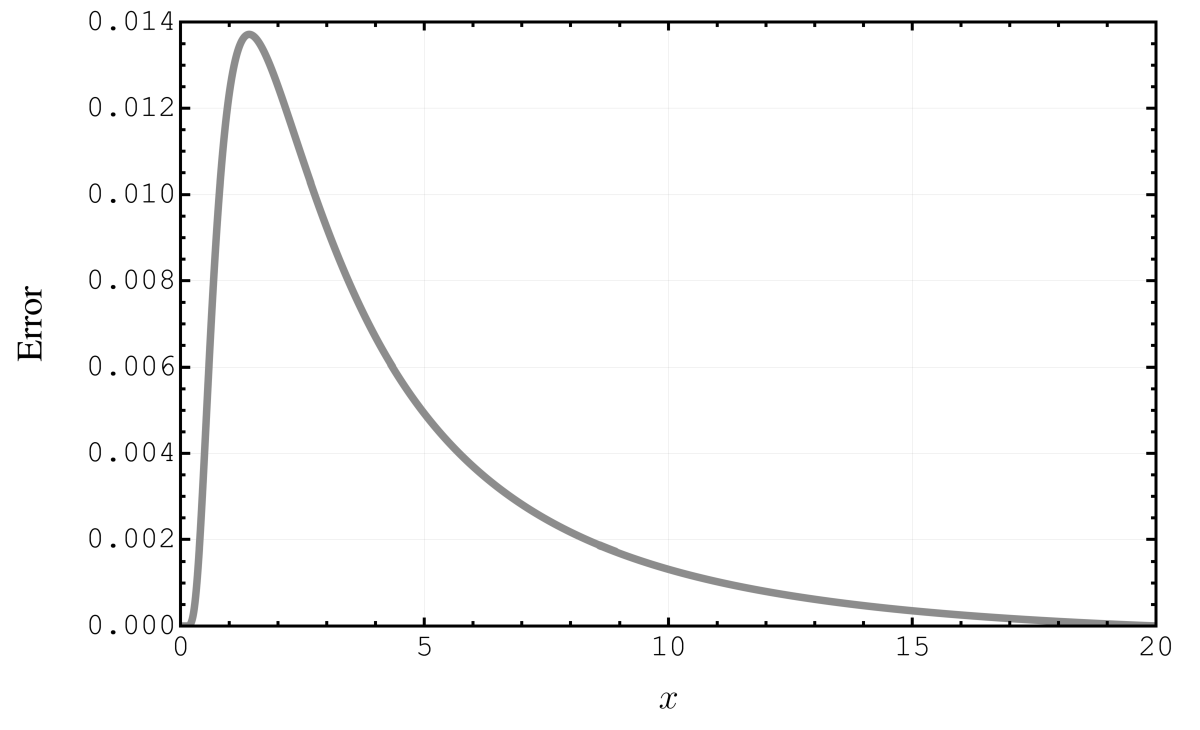}
        \caption{Corresponding lowerbound error.}
        \label{fig:QST-pdf-uprbnd5-err-A20}
    \end{subfigure}
    \caption{Quasi-stationary distribution's pdf, $q_{A}(x)$, its lowerbound $l_{A}^{(5)}(x)$, and the corresponding error---all as functions of $x\in[0,A]$ for $A=20$.}
    \label{fig:QST-pdf-lwrbnd5-perf-A20}
\end{figure}

\begin{figure}[h!]
    \centering
    \begin{subfigure}{0.48\textwidth}
        \centering
        \includegraphics[width=\linewidth]{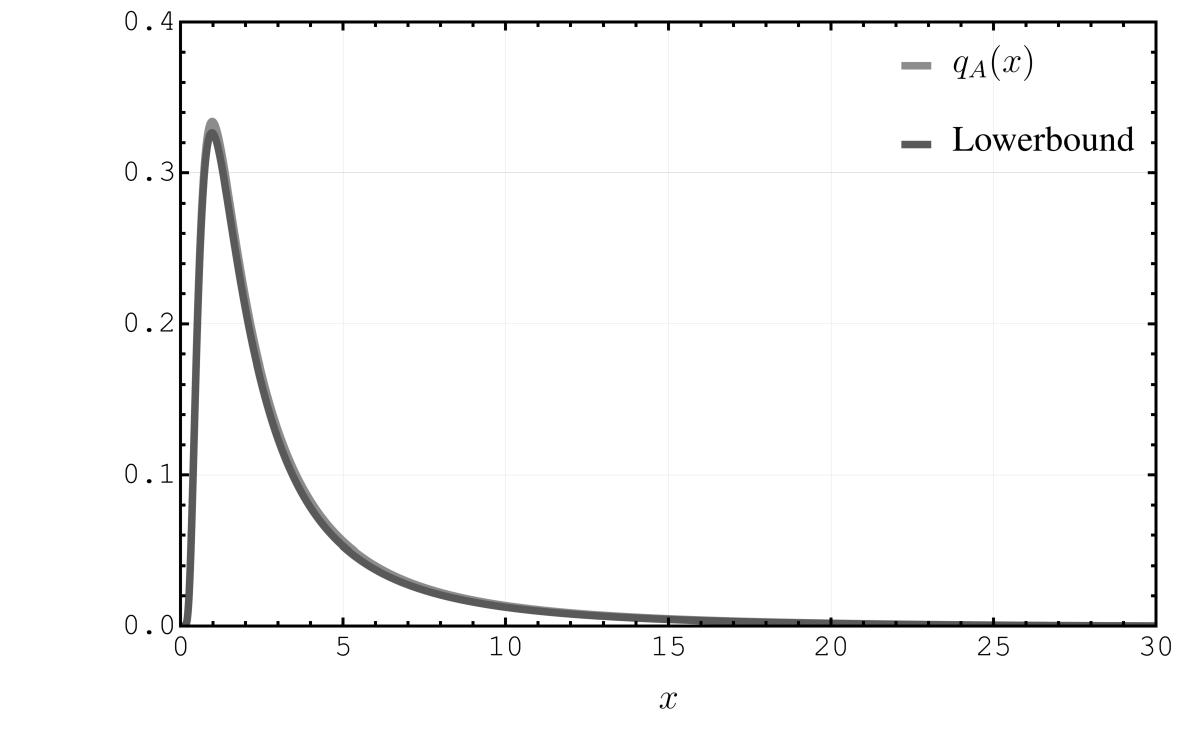}
        \caption{$q_{A}(x)$ and $l_{A}^{(5)}(x)$.}
        \label{fig:QST-pdf-lwrbnd5-A30}
    \end{subfigure}
    \hspace*{\fill}
    \begin{subfigure}{0.48\textwidth}
        \centering
        \includegraphics[width=\linewidth]{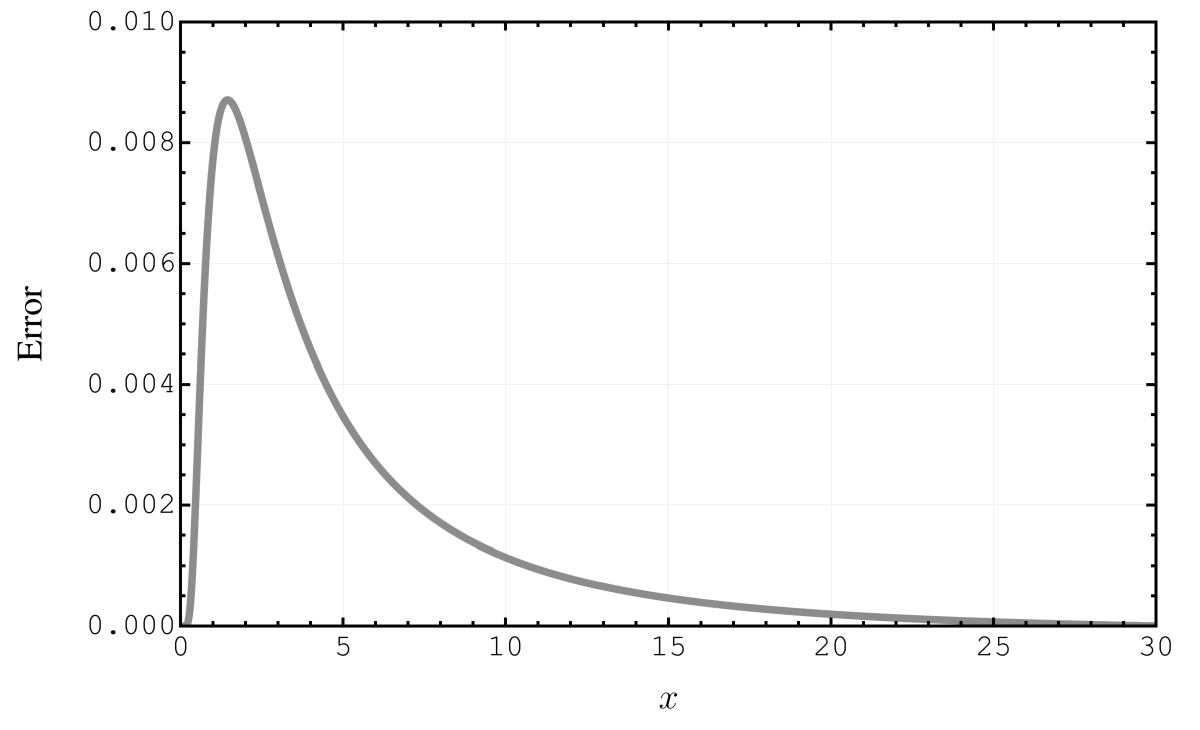}
        \caption{Corresponding lowerbound error.}
        \label{fig:QST-pdf-uprbnd5-err-A30}
    \end{subfigure}
    \caption{Quasi-stationary distribution's pdf, $q_{A}(x)$, its lowerbound $l_{A}^{(5)}(x)$, and the corresponding error---all as functions of $x\in[0,A]$ for $A=30$.}
    \label{fig:QST-pdf-lwrbnd5-perf-A30}
\end{figure}

Figures~\ref{fig:QST-pdf-lwrbnd5-perf-A20} and~\ref{fig:QST-pdf-lwrbnd5-perf-A30} show the performance of $l_{A}^{(5)}(x)$ as a function of $x\in[0,A]$ for $A=20$ and $A=30$. From the figures we see that the lowerbound is tighter than $l_{A}^{(1)}(x)$, for all $x\in[0,A]$.

This bound then gives a new upperbound for $Q_{A}(x)$, namely
\begin{align}\label{eq:QSD-cdf-uprbnd5-def}
Q_{A}(x)
&\le
e^{\tfrac{2}{A}}H(x)\left(\dfrac{A}{x}\right)^{\tfrac{2}{A}}
\eqqcolon U_{A}^{(5)}(x),
\;
x\in[0,A],
\;
A\ge\tilde{A},
\end{align}
which is tighter than $U_{A}^{(1)}(x)$ given by~\eqref{eq:QSD-cdf-uprbnd1-def}, because $\lambda>1/A$ for any $A>0$, by~\eqref{eq:lambda-dbl-ineq}. Moreover, unlike $U_{A}^{(1)}(x)$, the bound $U_{A}^{(5)}(x)$ is actually a cdf in itself: it is a smooth, strictly increasing function, and such that $0\le U_{A}^{(5)}(x)\le 1$ for all $x\in[0,A]$. The reason is the trivial inequality $\log(x)\le x-1$ valid for all $x>0$, so that
\begin{align*}
\log\left(\dfrac{A}{x}\right)
&<
\dfrac{A}{x}-1,
\;
\text{and}
\;
\dfrac{2}{A}\log\left(\dfrac{A}{x}\right)<\dfrac{2}{x}-\dfrac{2}{A}.
\end{align*}

\begin{figure}[h!]
    \centering
    \begin{subfigure}{0.48\textwidth}
        \centering
        \includegraphics[width=\linewidth]{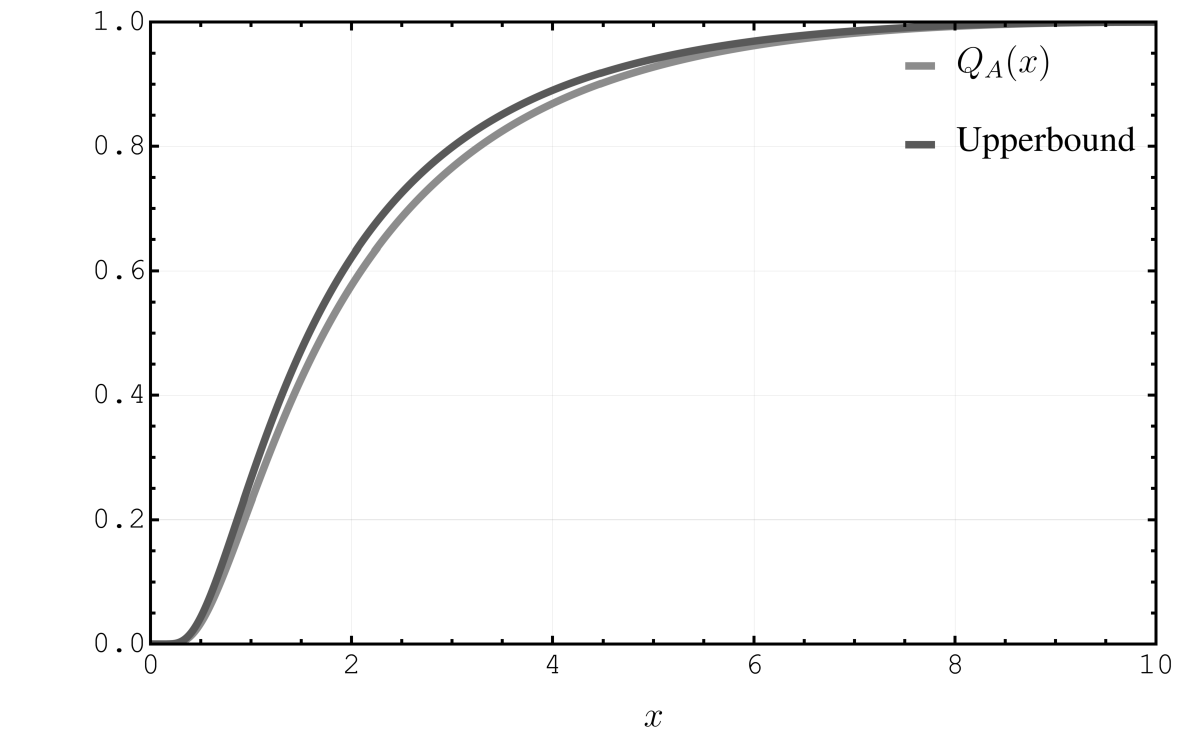}
        \caption{$Q_{A}(x)$ and $U_{A}^{(5)}(x)$.}
        \label{fig:QST-cdf-uprbnd2-A10}
    \end{subfigure}
    \hspace*{\fill}
    \begin{subfigure}{0.48\textwidth}
        \centering
        \includegraphics[width=\linewidth]{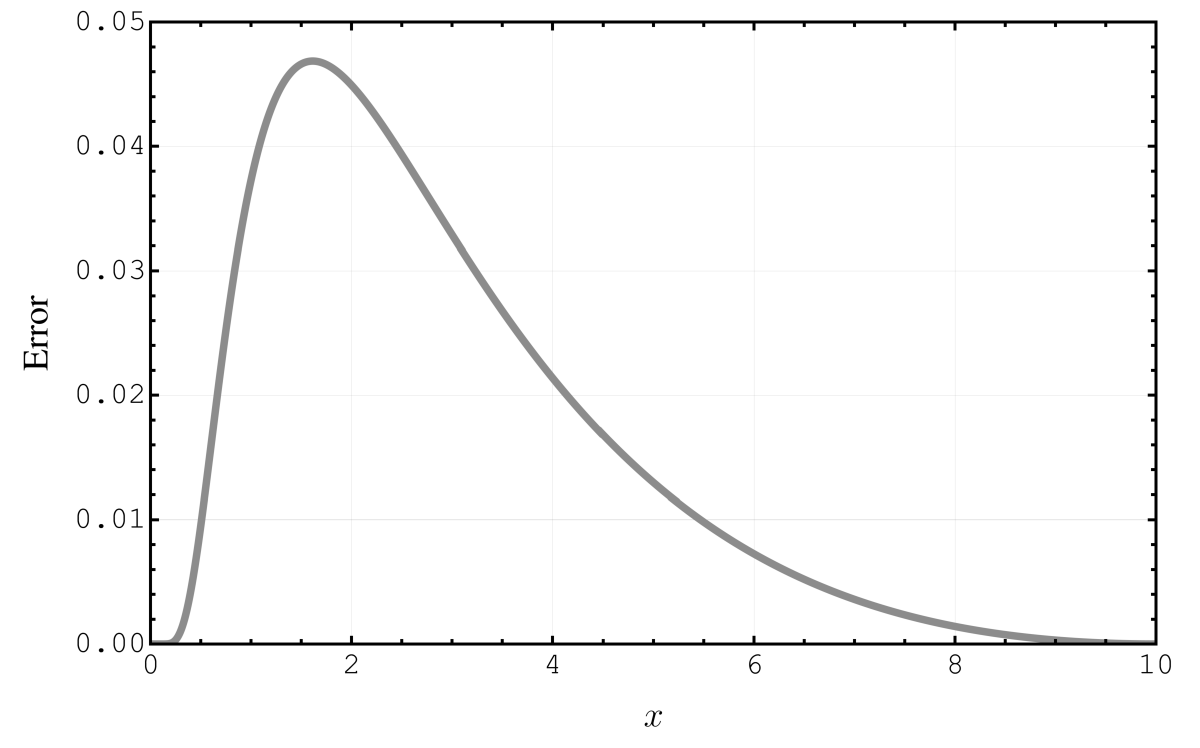}
        \caption{Corresponding upperbound error.}
        \label{fig:QST-cdf-uprbnd5-err-A10}
    \end{subfigure}
    \caption{Quasi-stationary distribution's cdf, $Q_{A}(x)$, its upperbound $U_{A}^{(5)}(x)$, and the corresponding error---all as functions of $x\in[0,A]$ for $A=10$.}
    \label{fig:QST-cdf-uprbnd5-perf-A10}
\end{figure}

\begin{figure}[h!]
    \centering
    \begin{subfigure}{0.48\textwidth}
        \centering
        \includegraphics[width=\linewidth]{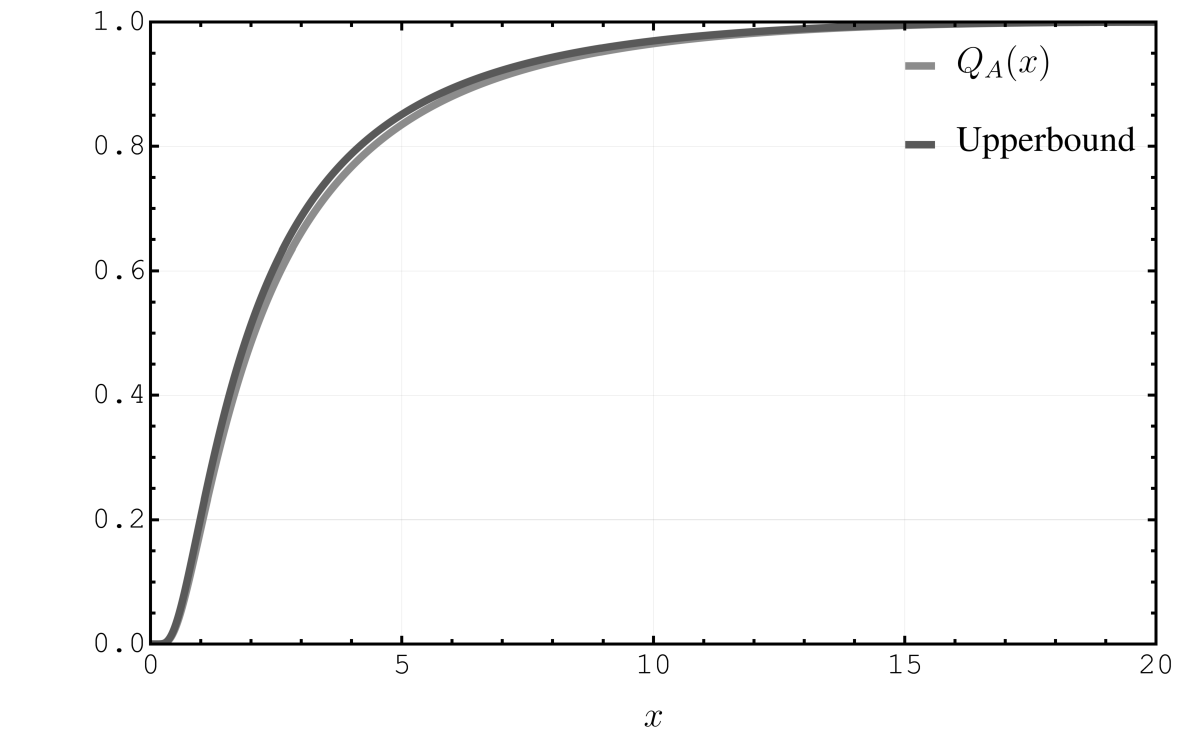}
        \caption{$Q_{A}(x)$ and $U_{A}^{(5)}(x)$.}
        \label{fig:QST-cdf-uprbnd5-A20}
    \end{subfigure}
    \hspace*{\fill}
    \begin{subfigure}{0.48\textwidth}
        \centering
        \includegraphics[width=\linewidth]{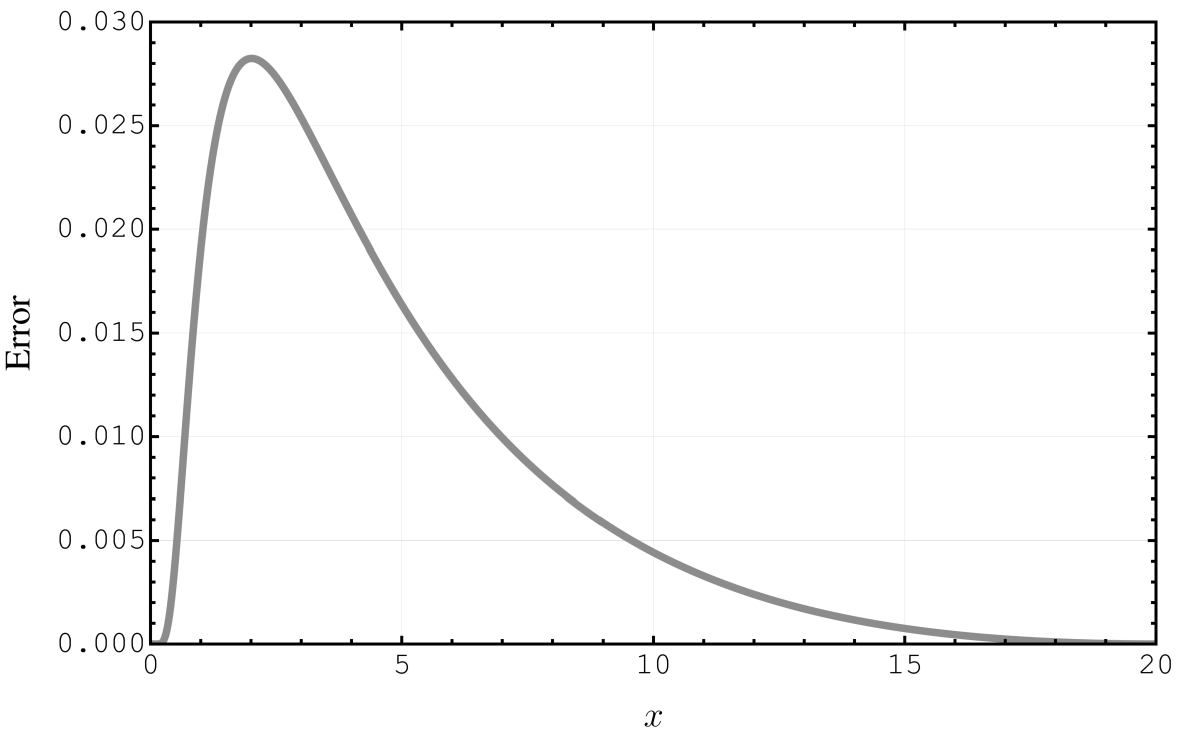}
        \caption{Corresponding upperbound error.}
        \label{fig:QST-cdf-uprbnd5-err-A20}
    \end{subfigure}
    \caption{Quasi-stationary distribution's cdf, $Q_{A}(x)$, its upperbound $U_{A}^{(5)}(x)$, and the corresponding error---all as functions of $x\in[0,A]$ for $A=20$.}
    \label{fig:QST-cdf-uprbnd5-perf-A20}
\end{figure}

\begin{figure}[h!]
    \centering
    \begin{subfigure}{0.48\textwidth}
        \centering
        \includegraphics[width=\linewidth]{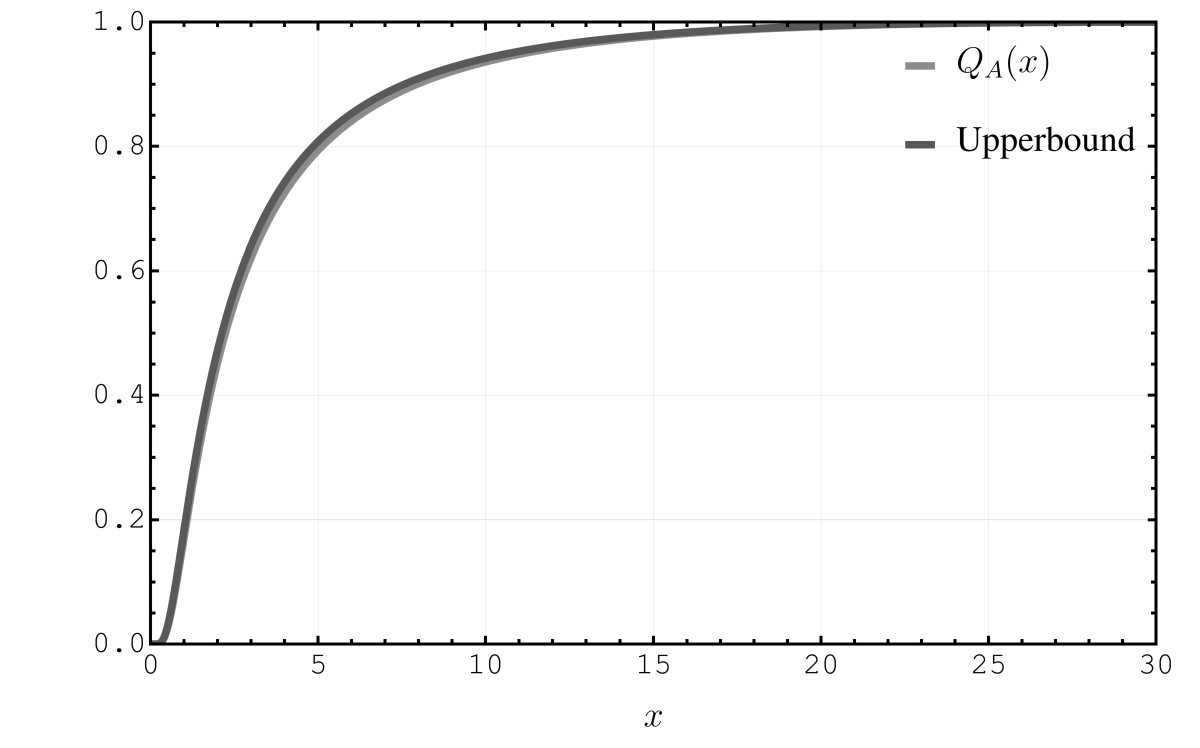}
        \caption{$Q_{A}(x)$ and $U_{A}^{(5)}(x)$.}
        \label{fig:QST-cdf-uprbnd5-A30}
    \end{subfigure}
    \hspace*{\fill}
    \begin{subfigure}{0.48\textwidth}
        \centering
        \includegraphics[width=\linewidth]{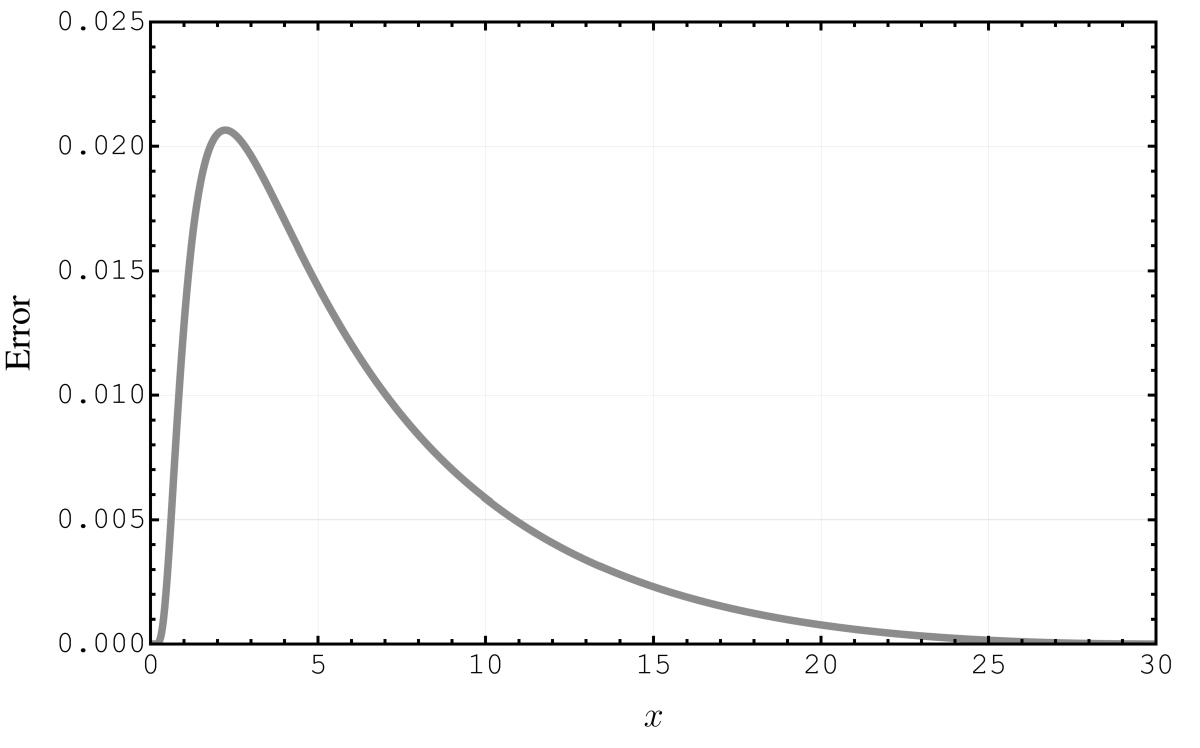}
        \caption{Corresponding upperbound error.}
        \label{fig:QST-cdf-uprbnd5-err-A30}
    \end{subfigure}
    \caption{Quasi-stationary distribution's cdf, $Q_{A}(x)$, its upperbound $U_{A}^{(5)}(x)$, and the corresponding error---all as functions of $x\in[0,A]$ for $A=30$.}
    \label{fig:QST-cdf-uprbnd5-perf-A30}
\end{figure}

Figures~\ref{fig:QST-cdf-uprbnd5-perf-A10},~\ref{fig:QST-cdf-uprbnd5-perf-A20}, and~\ref{fig:QST-cdf-uprbnd5-perf-A30} show the performance of $U_{A}^{(5)}(x)$ as a function of $x\in[0,A]$ for $A$ as low as $10$, $20$, and $30$. We note that $A=10$ is below $\tilde{A}$, and yet the bound $U_{A}^{(5)}(x)$ appears to work anyway, as can be seen in Figure~\ref{fig:QST-cdf-uprbnd5-perf-A10}. The figures show that the bound is extremely tight, for all $x\in[0,A]$, and gets even tighter as $A$ increases.

To draw a line under this section we point out that one can ``iterate'' either formula~\eqref{eq:QSD-pdf-qA-via-cdf-QA-v1} or formula~\eqref{eq:QSD-pdf-qA-via-cdf-QA-v2} to obtain new, even tighter, but more complicated, bounds for $q_{A}(x)$; these bounds, in turn, can then be integrated, of course, and thereby be ``converted'' into bounds for $Q_{A}(x)$. For example, from~\eqref{eq:QSD-cdf-lwrbnd1-def} plugged under the integral on the right of~\eqref{eq:QSD-pdf-qA-via-cdf-QA-v1}, we find
\begin{align*}
\begin{split}
\int_{0}^{x}Q_{A}(t)\,dt
&=
Q_{A}(x)\int_{0}^{x}\dfrac{Q_{A}(t)}{Q_{A}(x)}\,dt
\ge
Q_{A}(x) e^{\tfrac{2}{x}} \int_{0}^{x} e^{-\tfrac{2}{t}} dt
\\
&=
Q_{A}(x)\left\{x-2 e^{\tfrac{2}{x}}\E1\left(\dfrac{2}{x}\right)\right\},
\;
x\in[0,A],
\;
A>0,
\end{split}
\end{align*}
whence
\begin{align}\label{eq:QSD-pdf-qA-uprbnd6-def}
  q_A(x) & \le \dfrac{2}{x^2} Q_{A}(x)\left\{1-\lambda\left[x-2e^{\tfrac{2}{x}}\E1\left(\dfrac{2}{x}\right)\right]\right\}\indicator{x\in[0,A]}
  \eqqcolon u_{A}^{(6)}(x), \\
  &  \qquad x\in\mathbb{R}, \; A>0, \nonumber
\end{align}
and it is a very sharp upperbound for $q_{A}(x)$; recall that $\E1(x)$ here denotes the exponential integral function~\eqref{eq:Ei-func-def}. By virtue of integration, it can be translated into the lowerbound for $Q_{A}(x)$. Specifically, we find:
\begin{align*}
  L_A^{(6)}(x)
  & \coloneqq
    e^{\tfrac{2}{A}} H(x) \exp\left\{2\lambda\left[e^{\tfrac{2}{A}}\E1\left(\dfrac{2}{A}\right)-e^{\tfrac{2}{x}}\E1\left(\dfrac{2}{x}\right)\right]\right\}
\le
Q_{A}(x), \\
  & \qquad x\in[0,A], \; A>0,
\end{align*}
which is also a very tight lowerbound for $Q_{A}(x)$, clearly tighter than the lowerbound given by~\eqref{eq:QSD-cdf-lwrbnd1-def}. It is also worth pointing out that $L_{A}^{(6)}(x)$ is a cdf in itself, with the corresponding density's support being the interval $[0,A]$.

Similarly, we have
\begin{align*}
\dfrac{Q_{A}(x_1)}{Q_{A}(x_2)}
&\le
e^{\tfrac{2}{x_2}}e^{-\tfrac{2}{x_1}} x_{2}^{\tfrac{2}{A}} x_{1}^{-\tfrac{2}{A}},
\;
0\le x_1\le x_2\le A,
\;
A>0,
\end{align*}
and thus
\begin{align*}
\begin{split}
\int_{0}^{x}Q_{A}(t)\,dt
&=
Q_{A}(x)\int_{0}^{x}\dfrac{Q_{A}(t)}{Q_{A}(x)}\,dt
\le
Q_{A}(x) e^{\tfrac{2}{x}} x^{\tfrac{2}{A}} \int_{0}^{x} e^{-\tfrac{2}{t}} t^{-\tfrac{2}{A}} dt
\\
&=
2Q_{A}(x) e^{\tfrac{2}{x}} \left(\dfrac{x}{2}\right)^{\tfrac{2}{A}} \Gamma\left(-1+\dfrac{2}{A},\dfrac{2}{x}\right),
\;
x\in[0,A],
\;
A>0,
\end{split}
\end{align*}
so that
\begin{align*}
l_{A}^{(6)}(x)
&\coloneqq
\dfrac{2}{x^2} Q_{A}(x)\left\{1-2\lambda e^{\tfrac{2}{x}} \left(\dfrac{x}{2}\right)^{\tfrac{2}{A}} \Gamma\left(-1+\dfrac{2}{A},\dfrac{2}{x}\right)\right\}\indicator{x\in[0,A]}
\le
q_{A}(x), \\
& \qquad x\in\mathbb{R}, \; A>0,
\end{align*}
which is a very sharp lowerbound for $q_{A}(x)$; its sharpness though is offset by the presence of a special function, namely the incomplete Gamma function.

Translating the lowerbound $l_{A}^{(6)}(x)$ into an upperbound for $Q_{A}(x)$ is somewhat problematic, because $l_{A}^{(6)}(x)$ is difficult to integrate with respect to $x$, for
\begin{align*}
\begin{split}
-\log[Q_{A}(x)]
&\ge
\dfrac{2}{x}-\dfrac{2}{A}-4\lambda \int_{x}^{A} e^{\tfrac{2}{u}} \left(\dfrac{u}{2}\right)^{\tfrac{2}{A}} \Gamma\left(-1+\dfrac{2}{A},\dfrac{2}{u}\right) \dfrac{du}{u^2}\\
&=
\dfrac{2}{x}-\dfrac{2}{A}-2\lambda \int_{\tfrac{2}{A}}^{\tfrac{2}{x}} e^{t}\,t^{-\tfrac{2}{A}} \Gamma\left(-1+\dfrac{2}{A},t\right) dt.
\end{split}
\end{align*}

One way to handle the integral of the incomplete Gamma function is to use the integral representation
\begin{align*}
e^{t}\, t^{-a} \Gamma(a-1,t)
&=
\dfrac{1}{t}\int_{0}^{\infty} \dfrac{e^{-t y}}{(1+y)^{2-a}} dy,
\end{align*}
so that we readily get
\begin{align*}
\int_{\tfrac{2}{A}}^{\tfrac{2}{x}} e^{t}\,t^{-\tfrac{2}{A}} \Gamma\left(-1+\dfrac{2}{A},t\right) dt
&=
\int_{0}^{\infty}\left(\int_{\tfrac{2}{A}}^{\tfrac{2}{x}}\dfrac{e^{-ty}}{t}dt\right) \dfrac{dy}{(1+y)^{2-\tfrac{2}{A}}}
 \\
 & = 
\int_{0}^{\infty}\left[\E1\left(\dfrac{2}{A}y\right)-\E1\left(\dfrac{2}{x}y\right)\right]\dfrac{dy}{(1+y)^{2-\tfrac{2}{A}}}.
\end{align*}

Another option is to use the integral representation
\begin{align*}
e^{t}\, t^{-a} \Gamma(a-1,t)
&=
\dfrac{1}{t\Gamma(2-a)}\int_0^{\infty}\dfrac{y^{1-a} e^{-y}}{y+t}dy,
\end{align*}
whence
\begin{align*}
\begin{split}
\int_{\tfrac{2}{A}}^{\tfrac{2}{x}} e^{t}\, t^{-\tfrac{2}{A}} \Gamma\left(-1+\dfrac{2}{A},t\right) dt
&=
\int_{0}^{\infty}\left(\int_{\tfrac{2}{A}}^{\tfrac{2}{x}}\dfrac{dt}{t(t+y)}\right) y^{1-\tfrac{2}{A}} e^{-y} dy
\\
&=
\int_{0}^{\infty}\log\left(\dfrac{1+\tfrac{A}{2}y}{1+\tfrac{x}{2}y}\right)e^{-y} \dfrac{dy}{y^{\tfrac{2}{A}}}.
\end{split}
\end{align*}

We therefore arrive at the inequality
\begin{align*}
Q_{A}(x)
&\le
e^{\tfrac{2}{A}} H(x) \exp\left\{2\lambda \int_{0}^{\infty}\left[\E1\left(\dfrac{2}{A}y\right)-\E1\left(\dfrac{2}{x}y\right)\right]\dfrac{dy}{(1+y)^{2-\tfrac{2}{A}}} \right\}
\eqqcolon U_{A}^{(6)}(x), \\
  & \qquad x\in[0,A], \; A\ge\tilde{A},
\end{align*}
or equivalently
\begin{align*}
Q_{A}(x)
&\le
e^{\tfrac{2}{A}} H(x) \exp\left\{2\lambda \int_{0}^{\infty}\log\left(\dfrac{1+\tfrac{A}{2}y}{1+\tfrac{x}{2}y}\right)e^{-y} \dfrac{dy}{y^{\tfrac{2}{A}}} \right\}
\eqqcolon U_{A}^{(6)}(x), \\
  & \qquad x\in[0,A], \; A\ge\tilde{A},
\end{align*}
and $U_{A}^{(6)}(x)$ is a very sharp upperbound, but is not as simple as~\eqref{eq:QSD-cdf-uprbnd5-def}.

We presented quite a few lower- and upper-bounds for $q_A(x)$ as well as for $Q_A(x)$. The bounds vary in their complexity and tightness (accuracy): those that are more complex are more tight, and those that are simpler are less accurate. A good compromise, in our opinion, is bounds $l^{(5)}(x)$ and $u^{(1)}(x)$ for $q_A(x)$ and bounds $L^{(1)}(x)$ and $U^{(5)}(x)$ for $Q_A(x)$.

%

\section{Discussion}
\label{sec:discussion}
We now illustrate a few applications of the bounds obtained in the previous section for the quasi-stationary pdf $q_{A}(x)$ and cdf $Q_{A}(x)$.

Let us first try to use our bounds to quantify the difference between $Q_{A}(x)$ and $H(x)$ for $x\in(0,A)$ with $A\ge\tilde{A}$; otherwise, for $x\le 0$ we have $Q_{A}(x)=H(x)=0$ for any $A>0$, and for $x\ge A$ we have $Q_{A}(x)=1$ for any $A>0$, so $Q_{A}(x)-H(x)=1-e^{-\tfrac{2}{x}}\le 1-e^{-\tfrac{2}{A}}$. 

We start with the observation that
\begin{align*}
Q_{A}(x)-H(x)
&\le
e^{\tfrac{2}{A}} H(x) \left(\dfrac{A}{x}\right)^{\tfrac{2}{A}} - H(x)
\eqqcolon g(x;A),
\;
x\in\mathbb{R},
\;
A\ge\tilde{A},
\end{align*}
which is an immediate consequence of~\eqref{eq:QSD-cdf-uprbnd5-def}. Fix $A\ge\tilde{A}$, restrict $x$ to the interval $(0,A)$, and consider the function $g(x;A)$,
which is obviously a bounded and smooth function on $x\in(0,A)$. The first derivative of $g(x;A)$ with respect to $x$ is
\begin{align*}
\dfrac{\partial}{\partial x}
g(x;A)
&=
\dfrac{2}{x^2} g(x;A)-\dfrac{2}{A x}e^{\tfrac{2}{A}} H(x) \left(\dfrac{A}{x}\right)^{\tfrac{2}{A}}
,
\end{align*}
and it vanishes at $x_{A}^{*}\in(0,A)$ that solves the equation
\begin{align}\label{eq:Q-H-max-eqn1}
e^{\tfrac{2}{A}}\left(\dfrac{A}{x_{A}^{*}}\right)^{\tfrac{2}{A}}\left(1-\dfrac{x_{A}^{*}}{A}\right)
&=
1.
\end{align}

Equation~\eqref{eq:Q-H-max-eqn1} clearly has a unique solution, contained in the interval between 0 and $A$. To see this, fix $A$ and vary $x$ from $0$ up through $A$, and note that, on the on hand, the function $x\mapsto e^{2/A}(A/x)^{2/A}$ is strictly decreasing from $+\infty$ down to $e^{2/A}$, while, on the other hand, the function $x\mapsto 1/(1-x/A)$ is strictly increasing from 1 up to $+\infty$. The two functions definitely ``meet'' exactly once, at some point $x_{A}^{*}$ inside the interval $(0,A)$. Moreover, it is also clear that the function $g(x;A)$ is maximized at $x_{A}^{*}\in(0,A)$.

The exact solution $x_{A}^{*}$ to equation~\eqref{eq:Q-H-max-eqn1} is not possible to find analytically, due to the transcendental nature of the equation. However, it is not too difficult to upperbound $x_{A}^{*}$ via an elementary function of $A$. Specifically, since equation~\eqref{eq:Q-H-max-eqn1} is equivalent to the equation
\begin{align*} 
\dfrac{2}{A}+\dfrac{2}{A}\log(A)-\dfrac{2}{A}\log(x_{A}^{*})+\log\left(1-\dfrac{x_{A}^{*}}{A}\right)
&=
0,
\end{align*}
and $1-1/x\le\log(x)\le x-1$ for $x>0$, and $\log(1-x)\le -x$ for $x\in(0,1)$, we find
\begin{align*}
\log(x_{A}^{*})
&\ge
1-\dfrac{1}{x_{A}^{*}}
\;\;
\text{and}
\;\;
\log\left(1-\dfrac{x_{A}^{*}}{A}\right)
\le
-\dfrac{x_{A}^{*}}{A},
\end{align*}
so that
\begin{align*}
0
&\le
\dfrac{2}{A}+\dfrac{2}{A}\log(A)-\dfrac{2}{A}\left(1-\dfrac{1}{x_{A}^{*}}\right)-\dfrac{x_{A}^{*}}{A},
\end{align*}
whence
\begin{align*}
(x_{A}^{*})^2-2 x_{A}^{*} \log(A)-2
&\le
0,
\end{align*}
and thus
\begin{align*}
x_{A}^{*}
&\le
\log(A)+\sqrt{\log^2(A)+2}
\le
2\log(A)+\dfrac{1}{\log(A)},
\end{align*}
because $\sqrt{1-x}\le 1-x/2$ for $x\ge-1$.

Now, from~\eqref{eq:Q-H-max-eqn1} we see that
\begin{align*}
g(x_{A}^{*};A)
& =
e^{\tfrac{2}{A}} H(x_{A}^{*}) \left(\dfrac{A}{x_{A}^{*}}\right)^{\tfrac{2}{A}} - H(x_{A}^{*})
=
H(x_{A}^{*})\left(1-\dfrac{x_{A}^{*}}{A}\right)^{-1}-H(x_{A}^{*}) \\
& =
H(x_{A}^{*})\dfrac{x_{A}^{*}}{A-x_{A}^{*}}
\le
\dfrac{x_{A}^{*}}{A-x_{A}^{*}}.
\end{align*}

We can now claim that
\begin{align*}
\sup_{x\in(0,A)}[Q_{A}(x)-H(x)]
&\le
\dfrac{x_{A}^{*}}{A-x_{A}^{*}}
\le
\dfrac{2\log(A)+1/\log(A)}{A-2\log(A)-1/\log(A)}
=
O\left(\dfrac{\log(A)}{A}\right),
\end{align*}
which is consistent with~\eqref{eq:Q-to-H-unif-conv} proved by~\cite{Li+Polunchenko:SA2020}. By ``consistent'' we mean that~\cite{Li+Polunchenko:SA2020} merely showed the rate of convergence, i.e., the ``$\log(A)/A$'' under the ``big Oh'' symbol, but did not quantify the rate any more specifically. Here we fill in that void.


Next, let us try to use our bounds for $q_{A}(x)$ and $Q_{A}(x)$ to get a new lowerbound for $\lambda$. To that end, from~\eqref{eq:QA-int-lambda-link} and~\eqref{eq:QSD-cdf-uprbnd5-def} we have
\begin{align*}
\dfrac{1}{\lambda_{A}}
&=
\int_{0}^{A}Q_{A}(t)\,dt
\le
e^{\tfrac{2}{A}}\int_{0}^{A} H(x)\left(\dfrac{A}{x}\right)^{\tfrac{2}{A}}dx
=
2\left(\dfrac{A}{2}\right)^{\tfrac{2}{A}}e^{\tfrac{2}{A}}\Gamma\left(-1+\dfrac{2}{A},\dfrac{2}{A}\right),
\end{align*}
where $\Gamma(a,z)$ denotes the incomplete (upper) Gamma function. This gives at once
\begin{align*}
\lambda
&\ge
\dfrac{1}{2}\left(\dfrac{2}{A}\right)^{\tfrac{2}{A}}e^{-\tfrac{2}{A}}\left\{\Gamma\left(-1+\dfrac{2}{A},\dfrac{2}{A}\right)\right\}^{-1},
\end{align*}
which is a very sharp bound, but it has a special function involved in it.

To upperbound $\Gamma(a,z)$ we turn to~\cite[Proposition~2.11]{Pinelis:MIA2020}, i.e., the inequality
\begin{align*}
\Gamma(a,x)
&<
e^{-x}
\dfrac{x^a(1-a+x)}{(x-a)^2-a+2x},
\;
a<0,
\;
x>0,
\end{align*}
which gives
\begin{align*}
\Gamma\left(-1+\dfrac{2}{A},\dfrac{2}{A}\right)
&\le
\dfrac{A}{A+1}\left(\dfrac{A}{2}\right)^{\tfrac{2}{A}-1} e^{-\tfrac{2}{A}},
\end{align*}
and thus we arrive at the lowerbound
\begin{align*}
\dfrac{1}{A}+\dfrac{1}{A^2}
&\le
\lambda,
\end{align*}
which is clearly tighter than the left half of~\eqref{eq:lambda-dbl-ineq}.

From the lowerbound~\eqref{eq:QSD-pdf-qA-lwrbnd5-def} we find
\begin{align*}
\dfrac{x^2}{2}q_{A}(x)
&\ge
Q_{A}(x)\left(1-\dfrac{x}{A}\right)
,
\;
x\in[0,A],
\;
A\ge\tilde{A},
\end{align*}
so that
\begin{align*}
\dfrac{1}{2}\int_{0}^{A}x^2 q_{A}(x)dx
&\ge
\int_{0}^{A}Q_{A}(x)dx-\dfrac{1}{A}\int_{0}^{A}xQ_{A}(x)dx \\
& =
\int_{0}^{A}Q_{A}(x)dx-\dfrac{1}{A}\left(\dfrac{A^2}{2}-\dfrac{1}{2}\int_{0}^{A}x^2q_{A}(x)dx\right),
\end{align*}
whence
\begin{align*}
\int_{0}^{A}x^2 q_{A}(x)dx
&\ge
\dfrac{A(2-\lambda A)}{\lambda(A-1)},
\;
A\ge\tilde{A}.
\end{align*}

Now, from~\cite{Polunchenko:SA2017a} we have
\begin{align*}
\int_{0}^{A}x^2 q_{A}(x)dx
&=
\dfrac{\lambda^2 A^2-2\lambda A + 2}{\lambda+\lambda^2},
\end{align*}
so that we finally arrive at $\lambda^2 A^3-\lambda A^2-2\ge0$, whence
\begin{align*}
\lambda
&\ge
\dfrac{1+\sqrt{1+8/A}}{2A}
\ge
\dfrac{1}{A}+\dfrac{2}{A^2}\left(1-\dfrac{2}{A}\right),
\;
A\ge\tilde{A},
\end{align*}
which is tighter than our earlier result because $A\ge\tilde{A}$.

Finally, from the upperbound~\eqref{eq:QSD-pdf-qA-uprbnd6-def} and the fact that $q_{A}(A)=0$ we obtain
\begin{align*}
\lambda
&\le
\dfrac{1}{A}\left\{1-\dfrac{A}{2} e^{\tfrac{2}{A}}\E1\left(\dfrac{2}{A}\right)\right\}^{-1},
\end{align*}
which is a reasonably tight upperbound for $\lambda$.


\subsection{Application to Quickest Change-Point Detection}
Another important application of our bounds is new bounds for the delay exhibited by the Randomized Shiryaev--Roberts--Pollak (SRP) change-point detection method. Recall that the latter is identified with the stopping time $\mathcal{S}_{A}^{Q}$ given by~\eqref{eq:T-SRP-def}.

The worst-possible average detection delay exhibited by $\mathcal{S}_{A}^{Q}$ is given by
\begin{align*}
\SADD(\mathcal{S}_{A}^{Q})
&=
2\left\{e^{\tfrac{2}{A}}\E1\left(\dfrac{2}{A}\right)-\int_{0}^{A}e^{\tfrac{2}{t}}\E1\left(\dfrac{2}{t}\right)q_{A}(t)\,dt\right\},
\;
A>0;
\end{align*}
cf.~\cite{Feinberg+Shiryaev:SD2006}. It was shown by~\cite{Polunchenko:TPA2017} that the delay can also be expressed as
\begin{align*}
\SADD(\mathcal{S}_{A}^{Q})
&=
2\left\{e^{\tfrac{2}{A}}\E1\left(\dfrac{2}{A}\right)-1+2\lambda\int_{0}^{A}e^{\tfrac{2}{t}}\E1\left(\dfrac{2}{t}\right)Q_{A}(t)\dfrac{dt}{t}\right\},
\;
A>0,
\end{align*}
where $Q_{A}(x)$ is as in~\eqref{eq:QSD-cdf-answer-W0} or as in~\eqref{eq:QSD-cdf-answer-K}, and $\lambda$ is determined by~\eqref{eq:lambda-eqn}.

From~\eqref{eq:QSD-cdf-uprbnd5-def} we find
\begin{align*}
\int_{0}^{A}e^{\tfrac{2}{t}}\E1\left(\dfrac{2}{t}\right)Q_{A}(t)\dfrac{dt}{t}
&\le
e^{\tfrac{2}{A}} \left(\dfrac{A}{2}\right)^{\tfrac{2}{A}} \int_{0}^{A}\E1\left(\dfrac{2}{t}\right)\left(\dfrac{2}{t}\right)^{\tfrac{2}{A}}\dfrac{dt}{t} \\
& =
-\dfrac{A}{2}e^{\tfrac{2}{A}} \E1\left(\dfrac{2}{A}\right)+e^{\tfrac{2}{A}}\left(\dfrac{A}{2}\right)^{\tfrac{2}{A}+1}\Gamma\left(\dfrac{2}{A},\dfrac{2}{A}\right),
\end{align*}
because, on account of~\eqref{eq:E1-func-def}, we have
\begin{align*}
\dfrac{d}{dx}\left[\dfrac{x^{\varkappa+1}}{\varkappa+1}\E1(x)\right]
&=
x^{\varkappa}\E1(x)-\dfrac{x^{\varkappa}}{\varkappa+1} e^{-x},
\;
\varkappa\neq -1,
\end{align*}
and therefore
\begin{align*}
\int_{0}^{A}\E1\left(\dfrac{2}{t}\right)\left(\dfrac{2}{t}\right)^{\tfrac{2}{A}}\dfrac{dt}{t}
&=
\int_{\tfrac{2}{A}}^{+\infty} u^{\tfrac{2}{A}-1} \E1(u)\,du \\
& =
-\left(\dfrac{2}{A}\right)^{\tfrac{2}{A}-1}\E1\left(\dfrac{2}{A}\right)+\dfrac{A}{2}\Gamma\left(\dfrac{2}{A},\dfrac{2}{A}\right),
\end{align*}
for $\lim_{x\to0}[x^a \E1(x)]=0$, $a>0$, as given, e.g., by~\cite[Entry~3.5.3,~p.~193]{Geller+Ng:JRNBS1969}.

We thus have (for $A\ge\tilde{A}$)
\begin{align*}
\SADD(\mathcal{S}_{A}^{Q})
&\le
2\left\{\lambda A \left(\dfrac{A}{2}\right)^{\tfrac{2}{A}} e^{\tfrac{2}{A}}\Gamma\left(\dfrac{2}{A},\dfrac{2}{A}\right)-(\lambda A-1)e^{\tfrac{2}{A}} \E1\left(\dfrac{2}{A}\right)-1\right\},
\end{align*}
which is a very tight upperbound for $\SADD(\mathcal{S}_{A}^{Q})$, but it has the incomplete Gamma function in it. From the recurrence
\begin{align*}
\Gamma(a+1,x)
&=
a\Gamma(a,x)+x^{a}e^{-x}
\end{align*}
we find
\begin{align*}
\left(1-\dfrac{2}{A}\right)\Gamma\left(-1+\dfrac{2}{A},\dfrac{2}{A}\right)
&=
e^{-\tfrac{2}{A}}\left(\dfrac{2}{A}\right)^{-1+\tfrac{2}{A}}
-
\Gamma\left(\dfrac{2}{A},\dfrac{2}{A}\right),
\end{align*}
whence
\begin{align*}
e^{\tfrac{2}{A}}\left(\dfrac{A}{2}\right)^{\tfrac{2}{A}}\Gamma\left(-1+\dfrac{2}{A},\dfrac{2}{A}\right)
&\le
\dfrac{\lambda A-1}{2\lambda}+\dfrac{1}{\lambda A},
\end{align*}
and
\begin{align*}
\SADD(\mathcal{S}_{A}^{Q})
&\le
2(\lambda A-1)\left\{\dfrac{A}{2}-e^{\tfrac{2}{A}} \E1\left(\dfrac{2}{A}\right)\right\},
\;
A>\tilde{A},
\end{align*}
which is a reasonably tight bound, free of the incomplete Gamma function.



Now, from~\eqref{eq:QSD-pdf-qA-lwrbnd5-def} we find
\begin{align*}
\int_{0}^{A}e^{\tfrac{2}{t}}\E1\left(\dfrac{2}{t}\right)q_{A}(t)dt
&\ge
2\int_{0}^{A}e^{\tfrac{2}{t}}\E1\left(\dfrac{2}{t}\right)Q_{A}(t)\dfrac{dt}{t^2}-\dfrac{2}{A}\int_{0}^{A}e^{\tfrac{2}{t}}\E1\left(\dfrac{2}{t}\right)Q_{A}(t)\dfrac{dt}{t}
\end{align*}
whence
\begin{align*}
2\int_{0}^{A}e^{\tfrac{2}{t}}\E1\left(\dfrac{2}{t}\right)Q_{A}(t)\dfrac{dt}{t}
&\le
\dfrac{2\lambda A}{\lambda A - 1}
\left(1-2\int_{0}^{A}e^{\tfrac{2}{t}}\E1\left(\dfrac{2}{t}\right)Q_{A}(t)\dfrac{dt}{t^2} \right).
\end{align*}

On account of~\eqref{eq:QSD-cdf-lwrbnd1-def} we further obtain
\begin{align*}
\int_{0}^{A}e^{\tfrac{2}{t}}\E1\left(\dfrac{2}{t}\right)Q_{A}(t)\dfrac{dt}{t^2}
&\ge
e^{\tfrac{2}{A}}\int_{0}^{A}e^{\tfrac{2}{t}}\E1\left(\dfrac{2}{t}\right)H(t)\dfrac{dt}{t^2} \\
& =
e^{\tfrac{2}{A}}\int_{0}^{A}\E1\left(\dfrac{2}{t}\right)\dfrac{dt}{t^2}
=
\dfrac{1}{2}-\dfrac{1}{A}\,e^{\tfrac{2}{A}}\E1\left(\dfrac{2}{A}\right).
\end{align*}

Thus we arrive at
\begin{align*}
\SADD(\mathcal{S}_{A}^{Q})
&\le
2\left\{e^{\tfrac{2}{A}}\E1\left(\dfrac{2}{A}\right)-1+\dfrac{2\lambda}{\lambda A-1}e^{\tfrac{2}{A}}\E1\left(\dfrac{2}{A}\right)\right\},
\;
A\ge\tilde{A},
\end{align*}
which is a reasonably tight bound as well, but not as tight as the bound we obtained earlier.

\section{Concluding remarks} \label{sec:conclusion}

The primary goal of this work was to obtain new lower- and upper-bounds for the pdf, $q_{A}(x)$, as well as for the corresponding cdf, $Q_{A}(x)$, of the quasi-stationary distribution of the Generalized Shiryaev--Roberts (GSR) process. This is of importance in quickest change-point detection, which was the main motivation for us to initiate this work in the first place. We employed the latest monotonicity properties of the modified Bessel $K$ function and its functionals to achieve the goal, and offered a host of bounds---lower- and upper-bounds---for ratios $q_{A}(x)/h(x)$ and $Q_{A}(x)/H(x)$, where $h(x)$ and $H(x)$ are, respectively, the pdf and the cdf of the GSR process' stationary distribution; see~\eqref{eq:SR-StDist-def} and~\eqref{eq:SR-StDist-answer}.

However, there is more to the story: it is possible to quantify the difference between $q_{A}(x)$ and $h(x)$ as well as that between $Q_{A}(x)$ and $H(x)$ {\em directly}, i.e., obtain an exact closed-form expression for each. Specifically, the idea is to use the following identity
\begin{align*}
  \lefteqn{%
    \dfrac{\pi}{8}\dfrac{\Gamma(\sqrt{s}-1/2)}{\Gamma(1+2\sqrt{s})} W_{1,\sqrt{s}}(\max\{x,y\})M_{1,\sqrt{s}}(\min\{x,y\})} \\
  & = \dfrac{\pi}{2}\dfrac{xy}{4s-1}e^{-\tfrac{1}{2}(x+y)} +
  \int_0^{+\infty}\dfrac{\beta\sinh(\pi\beta)}{(1+4\beta^2)(s+\beta^2)}W_{1,\iu\beta}(x)W_{1,\iu\beta}(y)d\beta, \\
  & \qquad \text{valid for} \;\; s\in\mathbb{C}\setminus(-\infty,0)\cup\{1/4\} \;\; \text{and} \;\; x,y\ge0,
\end{align*}
where $M_{a,b}(z)$ is the Whittaker $M$ function, and $W_{a,b}(z)$ is the Whittaker $W$ function; cf.~\cite[Identity~(52),~p.~769]{Becker:JMP2004}. It should be noted that the case of $s=1/4$, while covered, has to be treated with care. On the one hand, the integral on the right {\em is} convergent for $s=1/4$; see~\cite{Becker:JMP2004} as well as our analysis below. However, on the other hand, as $s\to1/4$, two singularities emerge in the identity: one buried in the $\Gamma(\sqrt{s}-1/2)$ factor on the left, and the other one buried in the $1/(4s-1)$ factor on the right. The two are of the same order and, in the limit, as $s\to1/4$, end up ``undoing'' each other, yielding a finite expression for the limiting value of the integral.

With $\sqrt{s}=\xi/2$ and recalling that $\xi\in(0,1)$, the above identity takes the form
\begin{align*}
  \lefteqn{%
    \dfrac{\pi}{8}\dfrac{\Gamma(\xi/2-1/2)}{\Gamma(1+\xi)} W_{1,\tfrac{1}{2}\xi}\left(\dfrac{2}{x}\right)M_{1,\tfrac{1}{2}\xi}\left(\dfrac{2}{A}\right)} \\
  & = \dfrac{\pi}{2}\dfrac{2}{A}\dfrac{2}{x}\dfrac{1}{-8\lambda}e^{-\tfrac{1}{A}-\tfrac{1}{x}} \\
  & \quad + \int_0^{+\infty}\dfrac{4\beta\sinh(\pi\beta)}{(1+4\beta^2)(1-8\lambda+4\beta^2)}W_{1,\iu\beta}\left(\dfrac{2}{x}\right)W_{1,\iu\beta}\left(\dfrac{2}{A}\right)d\beta,
\\[1ex]
& \qquad x\in(0,A), \; A>\tilde{A}.
\end{align*}
To make use of this identity we need to recall
\begin{align*}
\lambda A\,\Gamma\left(\dfrac{\xi-1}{2}\right)W_{0,\tfrac{1}{2}\xi}\left(\dfrac{2}{A}\right)\,M_{1,\tfrac{1}{2}\xi}\left(\dfrac{2}{A}\right)
=
-\Gamma(\xi+1),
\end{align*}
shown in~\cite{Polunchenko+Pepelyshev:SP2018}.

From this and~\eqref{eq:SR-StDist-answer}, we arrive at
\begin{align}\label{eq:q-h-diff}
  \lefteqn{q_{A}(x)-h(x)} \\
 & =
-
\dfrac{8}{\pi}\dfrac{A\lambda}{x}e^{\tfrac{1}{A}-\tfrac{1}{x}}\int_{0}^{+\infty}\dfrac{4\beta\sinh(\pi\beta)}{(1+4\beta^2)(1-8\lambda+4\beta^2)}W_{1,\iu\beta}\left(\dfrac{2}{x}\right)W_{1,\iu\beta}\left(\dfrac{2}{A}\right)d\beta, \nonumber
\end{align}
for any $x\in[0,A]$ and $A>\tilde{A}$.

We have thus obtained an {\em exact} representation of the discrepancy between $q_{A}(x)$ and $h(x)$ valid for any fixed $x\in[0,A]$ and $A>\tilde{A}$. The integral on the right is absolutely convergent for any fixed $x\in(0,A]$. This is mentioned in~\cite{Becker:JMP2004}.

By integrating both sides of~\eqref{eq:q-h-diff} with respect to $x$ from 0 up through $t\in(0,A]$ we obtain the following:
\begin{align}\label{eq:Q-H-diff}
 \lefteqn{(0\le)\;Q_{A}(x)-H(x)} \\
  & =
-\dfrac{8}{\pi}A\lambda e^{\tfrac{1}{A}-\tfrac{1}{x}}\int_{0}^{+\infty}\dfrac{4\beta\sinh(\pi\beta)}{(1+4\beta^2)(1-8\lambda+4\beta^2)}W_{0,\iu\beta}\left(\dfrac{2}{x}\right)W_{1,\iu\beta}\left(\dfrac{2}{A}\right)d\beta, \nonumber
\end{align}
where we used the property of the Whittaker $W$ function that
$$
\dfrac{d}{dx}\left(e^{-x/2}W_{0,b}(x)\right)=-\dfrac{e^{-x/2}}{x}W_{1,b}(x).
$$

Formulae~\eqref{eq:q-h-diff} and~\eqref{eq:Q-H-diff} may lead to new bounds for $|q_{A}(x)-h(x)|$ and for $(0\le)\;Q_{A}(x)-H(x)=|Q_{A}(x)-H(x)|$. Getting these bounds comes down to whether or not we can tightly upperbound the integral on the right of~\eqref{eq:q-h-diff} and that on the right of~\eqref{eq:Q-H-diff}. To that end, the problem is that both integrals involve the Whittaker $W$ function, which is a somewhat exotic function in the sense that its properties are not as well understood as those of the modified Bessel $K$ function.
Hence, using~\eqref{eq:q-h-diff} and~\eqref{eq:Q-H-diff} to get new bounds for $|q_{A}(x)-h(x)|$ and for $(0\le)\;Q_{A}(x)-H(x)=|Q_{A}(x)-H(x)|$ will likely prove to be a challenge. Nonetheless, it does seem worthwhile, and, if and once solved, the results will be published in a separate paper.


\section*{Acknowledgement}
We thank Prof. Sven Knoth of the Helmut Schmidt University for his attention to the paper.

The effort of A.\,S.~Polunchenko was partially supported by the Simons Foundation via a Collaboration Grant in Mathematics under Award \#\,304574.

\bibliography{main,physics,special-functions,finance,stochastic-processes,differential-equations}


\end{document}